\documentclass[12pt,twoside]{article}
\usepackage{ambp}
\usepackage{amsfonts}

\usepackage[french]{babel}
\usepackage{epsfig}
\usepackage[changebar]{}
\usepackage[latin1]{inputenc}
\usepackage{makeidx}
\usepackage[french]{babel}

\begin{document}

\JournalVolume{10}
\JournalYear{2002}


\font \sevenrm=cmr7 \font \eightrm=cmr8 \font \bigrm=cmr10 scaled
\magstep1 \font \Bigrm=cmr10 scaled \magstep2 \font \biggrm=cmr10
scaled \magstep3 \font \Biggrm=cmr10 scaled \magstep4
\font \eightbf=cmbx8 \font \bigbf=cmbx10 scaled \magstep1 \font
\Bigbf=cmbx10 scaled \magstep2 \font \biggbf=cmbx10 scaled
\magstep3 \font \Biggbf=cmbx10 scaled \magstep4
%
\font \tengoth=eufm10 \font\eightgoth=eufm8 \font \sevengoth=eufm7
\font \fivegoth=eufm5 \font \biggoth=eufm10 scaled \magstep1 \font
\Biggoth=eufm10 scaled \magstep2
\newfam\gothfam
\textfont \gothfam=\tengoth \scriptfont \gothfam=\sevengoth
\scriptscriptfont \gothfam=\fivegoth
\def\goth{\fam\gothfam\tengoth}

%
\font \tenmath=msbm10 \font \sevenmath=msbm7 \font \fivemath=msbm5

\newfam\mathfam
\textfont \mathfam=\tenmath \scriptfont \mathfam=\sevenmath
\scriptscriptfont \mathfam=\fivemath
\def\math{\fam\mathfam\tenmath}
%
%
%
%
\def\ssq{\vskip 1mm}
\def\msq{\vskip 2mm}
\def\bsq{\vskip 3mm}
%
%
\def\uple#1#2{#1_1,\ldots ,{#1}_{#2}}
\def \restr#1{\mathstrut_{\textstyle |}\raise-6pt\hbox{$\scriptstyle #1$}}
\def \srestr#1{\mathstrut_{\scriptstyle |}\hbox to -1.5pt{}
\raise-4pt\hbox{$\scriptscriptstyle #1$}}
\def \inver{^{-1}}
\def\surj#1{\mathop{\hbox to #1 mm{\rightarrowfill\hskip 2pt\llap{$\rightarrow$}}}\limits}
\def\R{{\math R}}
\def\C{{\math C}}
\def\Q{{\math Q}}
\def\H{{\math H}}
\def\N{{\math N}}
\def\Z{{\math Z}}
\def\I{{\math I}}
\def\permuc#1#2#3#4{#1#2#3#4+#1#3#4#2+#1#4#2#3}
\def\fleche#1{\mathop{\hbox to #1 mm{\rightarrowfill}}\limits}
%
%
\def \g#1{\hbox{\tengoth #1}}
\def\mg#1{\hbox{\eightgoth #1}}
\def \sg#1{\hbox{\sevengoth #1}}
\def \ssg#1{\hbox{\fivegoth #1}}
\def\Cal #1{{\cal #1}}
%
%
\def \mop#1{\mathop{\hbox{\rm #1}}\nolimits}
\def \smop#1{\mathop{\hbox{\sevenrm #1}}\nolimits}
\def \ssmop#1{\mathop{\hbox{\fiverm #1}}\nolimits}
\def \mopl#1{\mathop{\hbox{\rm #1}}\limits}
\def \smopl#1{\mathop{\hbox{\sevenrm #1}}\limits}
\def \ssmopl#1{\mathop{\hbox{\fiverm #1}}\limits}
\let\wh=\widehat
\def\wwedge{\wedge\cdots\wedge}
\long\def\dessin#1#2{\null
               \bsq
               \begin{center}
                \epsfig{file=#2, height=#1, clip=}
                \end{center}
                \bsq
                \goodbreak}

\def\dm{\frac 12}
\def\hdm{\frac \hbar 2}
\let\wt=\widetilde
\let\bord=\partial

%
\begin{center}
{\Large \textbf{ Cohomologie tangente et cup-produit pour la
quantification de Kontsevich}\\\vspace{1cm}
Dominique Manchon\\
Charles Torossian}
\end{center}

\FirstAuthor{Dominique Manchon} \FirstDepartment{CNRS-UMR 6620}
\FirstInstitution{Universit\'e Blaise Pascal}
\FirstStreetAddress{24 avenue des Landais} \FirstCity{63177
Aubi\`ere Cedex} \FirstCountry{France}
\FirstEmail{manchon@math.univ-bpclermont.fr}
%
%
\SecondAuthor{Charles Torossian} \SecondDepartment{CNRS - UMR
8553}
\SecondInstitution{Ecole Normale Sup\'erieure}
\SecondStreetAddress{45 rue d'Ulm} \SecondCity{75230 Paris Cedex
05} \SecondCountry{France}
\SecondEmail{Charles.Torossian@ens.fr}
%
%
%
%


\ShortAuthor{D. Manchon, Ch. Torossian}

\ShortTitle{Cohomologie tangente et cup-produit}

%
%
%


\begin{abstract}
On a flat manifold $M=\R^d$, M. Kontsevich's formality
quasi-isomorphism is compatible with cup-products on tangent
cohomology spaces, in the sense that for any formal Poisson
$2$-tensor $\hbar\gamma$ the derivative at $\hbar\gamma$ of the
quasi-isomorphism induces an isomorphism of graded commutative
algebras from Poisson cohomology space to Hochschild cohomology
space relative to the deformed multiplication built from
$\hbar\gamma$ via the quasi-isomorphism. We give here a detailed
proof of this result, with signs and orientations precised.
\end{abstract}


\section{Introduction}

L'existence d'un \'etoile-produit \cite{BF} sur une vari\'et\'e de
Poisson quelconque est la cons\'e\-quence d'un r\'esultat plus
profond : le th\'eor\`eme de formalit\'e de M. Kontsevich \cite K,
qui \'enonce l'existence d'un quasi-isomorphisme $L_\infty$ entre
les deux alg\`ebres de Lie diff\'erentielles gradu\'ees
naturellement attach\'ees \`a une vari\'et\'e diff\'erentiable $M$~:
l'alg\`ebre de Lie diff\'eren\-tielle gradu\'ee $\g g_1$ des
multi-champs de vecteurs (ou tenseurs contravariants) sur $M$ (avec
diff\'erentielle nulle et crochet de Schouten), et l'alg\`ebre de
Lie diff\'erentielle gradu\'ee $\g g_2$ des op\'erateurs
multi-diff\'erentiels sur $M$ (avec diff\'erentielle de Hochschild
et crochet de Gerstenhaber). \ssq Un tel quasi-isomorphisme
$L_\infty$ fournit un proc\'ed\'e canonique et parfaitement
explicite pour produire un \'etoile-produit $*=*_{\gamma}$ \`a
partir d'un $2$-tenseur de Poisson formel $\gamma$ sur la
vari\'et\'e, et tout \'etoile-produit est \'equivalent \`a un
\'etoile-produit obtenu de cette fa\c con (\cite {K} \S ~4.4, \cite
{AMM} \S~A.2). \ssq Les graduations de $\g g_1$ et $\g g_2$ sont
telles qu'un \'el\'ement homog\`ene de degr\'e $n$ dans $\g g_1$
(resp. $\g g_2$) est un $(n+1)$-champ de vecteur (resp. un
op\'erateur $(n+1)$-diff\'erentiel). Nous consid\'erons comme dans
\cite K les espaces d\'ecal\'es $\g g_1[1]$ et $\g g_2[1]$ comme des
vari\'et\'es formelles gradu\'ees point\'ees, en ce sens que pour
$i=1,2$ la structure d'alg\`ebre de Lie diff\'erentielle gradu\'ee
d\'efinit une cod\'erivation $Q^i$ de degr\'e $1$ de la cog\`ebre
sans co-unit\'e $S^+(\g g_i[1])$ qui v\'erifie l'\'equation ma\^\i
tresse~:
\begin{equation}\label{maitresse}
[Q^i,Q^i]=0.
\end{equation}
L'espace d\'ecal\'e $\g g_i[1]$ d\'esigne l'espace $\g g_i$ dans
lequel la graduation a augment\'e d'une unit\'e~: ainsi un
\'el\'ement homog\`ene de degr\'e $n$ dans $\g g_1[1]$ est un
$(n+2)$-champ de vecteur, et un \'el\'ement homog\`ene de degr\'e
$n$ dans $\g g_2[1]$ est un op\'erateur $(n+2)$-diff\'erentiel.

\begin{theoreme}
Il existe un quasi-isomorphisme $L_\infty$ de la vari\'et\'e
formelle gradu\'ee point\'ee $\g g_1[1]$ vers la vari\'et\'e
formelle gradu\'ee point\'ee $\g g_2[1]$, c'est-\`a-dire un
morphisme de cog\`ebres~:
$$\Cal U:S^+(\g g_1[1])\fleche 8 S^+(\g g_2[1])$$
tel que~:
\begin{equation}\label{formalite}
\Cal U\circ Q^1=Q^2\circ U,
\end{equation}
et tel que la restriction $\Cal U_1$ de $\Cal U$ \`a $\g g_1[1]$ est
un quasi-isomorphisme de complexes de cocha\^\i
nes\footnote{D'apr\`es \cite {AMM} on doit remplacer le crochet de
Schouten par l'oppos\'e du crochet pris dans l'ordre inverse. C'est
ce que nous ferons. Ce crochet co\"\i ncide avec le crochet de
Schouten modulo un signe moins lorsque deux \'el\'ements impairs
sont en jeu (cf. \S\ 4.2).}.
\end{theoreme}

Le th\'eor\`eme de formalit\'e est reli\'e \`a la quantification par
d\'eformation de la mani\`ere sui\-vante~: par propri\'et\'e
universelle des cog\`ebres cocommutatives colibres, les
cod\'erivations $Q^i$ et le $L_\infty$-morphisme $\Cal U$ sont
enti\`erement d\'etermin\'es par leurs coefficients de Taylor~:
\begin{eqnarray}
Q_k^i:S^k(\g g_i[1])         &\fleche 8 \g g_i[2]    \\
                \Cal U_k:S^k(\g g_1[1]) &\fleche 8 \g g_2[1],
\end{eqnarray}
$k\ge 1, \, i=1,2,$ obtenus en composant $Q^i$ et $\Cal U$ \`a
droite avec la projection canonique~: $\pi : S^+(\g g_i)\surj 5 \g
g_i$ (resp. $\pi : S^+(\g g_2)\surj 5 \g g_2$). Soit $\g
m=\hbar\R[[\hbar]]$ la limite projective des alg\`ebres nilpotentes
de dimension finie $\g m_r=\hbar\R[[\hbar]]/\hbar^r\R[[\hbar]]$.
Soit
$\hbar\gamma=\hbar(\gamma_0+\hbar\gamma_1+\hbar^2\gamma_2+\cdots)$
un $2$-tenseur de Poisson formel infinit\'esimal, c'est-\`a-dire une
solution dans $\g g_1\widehat\otimes \g m$ de l'\'equation de
Maurer-Cartan~:
\begin{equation}\label{mc}
\hbar d\gamma-\frac 12[\hbar\gamma,\hbar\gamma]=0,
\end{equation}
qui s'\'ecrit aussi de mani\`ere plus g\'eom\'etrique~:
\begin{equation}\label{mcg}
Q^1(e^{.\hbar\gamma}-1)=0,
\end{equation}
o\`u $e^{.\hbar\gamma}-1$ est un \'el\'ement de type groupe dans la
$\g m$-cog\`ebre topologique $S^+(\g g_1[1])\widehat\otimes \g m$.
Alors $\Cal U(e^{.\hbar\gamma}-1)$ est de type groupe dans la
cog\`ebre topologique $S^+(\g g_2[1])\widehat\otimes \g m$, et
donc~:
\begin{equation}\label{tg}
\Cal U(e^{.\hbar\gamma}-1)=e^{.\hbar\tilde \gamma}-1
\end{equation}
avec~:
\begin{equation}\label{tg2}
\hbar\tilde \gamma=\sum\limits_{k\ge 1}{\hbar^k\over k!}\Cal
U_k(\gamma^{.k}).
\end{equation}
Comme $Q^2$ s'annule en $e^{.\hbar\tilde\gamma}-1$ l'\'el\'ement
$\hbar\tilde \gamma$ verifie l'\'equation de Maurer-Cartan dans $\g
g_2\widehat\otimes \g m$~:
\begin{equation}\label{mc2}
\hbar d\tilde\gamma -{1\over 2}[\hbar\tilde \gamma,\hbar\tilde \gamma]=0.
\end{equation}
Soit $m$ l'op\'erateur bidiff\'erentiel de
multiplication~:$f\otimes g\mapsto fg$, et soit $*=m+\hbar\tilde
\gamma$. On rappelle (\cite {AMM} \S\ IV.3) que le cobord de
Hochschild est donn\'e par $d=-[m,-]$. L'\'equation de
Maurer-Cartan pour $\hbar\wt\gamma$ est donc \'equivalente \`a
l'\'equation~:
\begin{equation}\label{assoc}
[*,*]=0,
\end{equation}
qui dit que $*$ est un produit associatif sur
$C^\infty(M)[[\hbar]]$. \ssq L'\'el\'ement de type groupe $e^{.\hbar
\gamma}-1$ d\'esigne ``le $\g m$-point $\hbar \gamma$ sur la
vari\'et\'e formelle $\g g_1[1]$''. Cette expression prend tout son
sens g\'eom\'etrique si on pense \`a $e^{.\hbar \gamma}-1$ comme \`a
la diff\'erence entre les deux mesures de Dirac
$\delta_{\hbar\gamma}-\delta_0$. Dire que le champ de vecteurs
impair $Q^1$ s'annule au point $\hbar\gamma$, c'est dire
pr\'ecis\'ement que la cod\'erivation $Q^1$ s'annule sur $e^{.\hbar
\gamma}-1$. \ssq Nous nous int\'eressons maintenant \`a la
diff\'erentielle $\Cal U'_{\hbar\gamma}$ du morphisme de
vari\'et\'es formelles $\Cal U$. L'espace tangent \`a $\g g_1[1]$ en
$\hbar \gamma$ s'identifie \`a $\g g_1[1]\widehat\otimes \g m$.  La
lin\'earisation du champ de vecteurs impair $Q^1$ qui s'annule au
point $\hbar\gamma$ donne un champ de vecteurs impair
$Q^{\hbar\gamma}$ sur cet espace tangent, donn\'e par~:
\begin{equation}\label{vect}
Q^{\hbar\gamma}(\hbar\delta)=\sum\limits_{n\ge 0}{\hbar^{n+1}\over
n!}Q^1_{n+1} (\delta.\gamma^{.n}),
\end{equation}
comme on peut le voir en
extrayant le terme lin\'eaire en $\delta$ dans l'expresssion
$Q^1(e^{.\hbar\gamma+\hbar\delta}-1)$. L'\'equation ma\^\i tresse~:
\begin{equation}\label{maitresse2}
[Q^1,Q^1](e^{.\hbar\gamma+\hbar\delta}-1)=0
\end{equation}
permet facilement de d\'eduire l'\'equation ma\^\i tresse pour le
champ lin\'earis\'e~:
\begin{equation}\label{maitresse3}
[Q^{\hbar\gamma},Q^{\hbar\gamma}]=0.
\end{equation}
Dans notre cas particulier o\`u tous les coefficients de Taylor de
$Q^1$ sont nuls sauf le deuxi\`eme on obtient~:
\begin{equation}\label{vect2}
Q^{\hbar\gamma}(\hbar\delta)=Q^1_2(\hbar\gamma.\hbar\delta)=\hbar^2[\delta,\gamma].
\end{equation}
(Voir \cite {AMM} \S\ II.4 et IV.1). On lin\'earise
de m\^eme le champ de vecteurs impair $Q^2$ qui s'annule en
$\hbar\tilde\gamma=\Cal U(\hbar\gamma)$. On obtient ainsi un champ
de vecteurs impair $Q^{\hbar\tilde\gamma}$ de carr\'e nul sur
l'espace tangent \`a $\g g_2[1]$ en $\hbar\tilde\gamma=\Cal
U(\hbar\gamma)$, qui s'\'ecrit~:
\begin{equation}\label{vect3}
Q^{\hbar\tilde\gamma}(\hbar\delta)=Q^2_1(\hbar\delta)+Q^2_2(\hbar\tilde\gamma.
\hbar\delta)=\hbar[\delta,*].
\end{equation}
Les deux espaces tangents ci-dessus sont donc ainsi munis d'une
structure de complexe de cocha\^\i nes. La d\'eriv\'ee $\Cal
U'_{\hbar\gamma}$ du $L_\infty$-morphisme $\Cal U$ est un
quasi-isomorphisme de complexes du premier espace tangent vers le
deuxi\`eme, et s'exprime en \'etendant pas $\g m$-lin\'earit\'e la
formule ~:
\begin{equation}\label{derivee}
\Cal U'_{\hbar\gamma}(\delta)=\sum\limits_{n\ge 0}{\hbar^n\over
n!}\Cal U_{n+1} (\delta.\gamma^{.n}),
\end{equation}
pour $\delta \in \g g_1[1]$. L'op\'erateur de cobord
$Q^{\hbar\gamma}=[-,\hbar\gamma]=-[\hbar\gamma, -]_{\hbox{\sevenrm
Shouten}}$ est une d\'erivation gradu\'ee pour le produit
ext\'erieur $\wedge$ des multi-champs de vecteurs \'etendu \`a $\g
g_1[1]\wh\otimes\g m$ par $\g m$-lin\'earit\'e, par d\'efinition
m\^eme du crochet de Schouten. Le produit ext\'erieur induit donc un
produit associatif et commutatif (au sens gradu\'e), que nous
noterons $\cup$, sur l'espace de cohomologie $H_{\hbar\gamma}$ du
premier espace tangent. \ssq On introduit un produit associatif
gradu\'e sur le deuxi\`eme espace tangent en posant pour chaque
op\'erateur $k_1$-diff\'erentiel $t_1$ et pour chaque op\'erateur
$k_2$-diff\'erentiel $t_2$~:
\begin{equation}\label{cup}
(t_1\cup t_2)(a_1\otimes\cdots\otimes a_{k_1+k_2})=t_1(a_1\otimes\cdots
\otimes a_{k_1})*t_2(a_{k_1+1}\otimes\cdots\otimes
a_{k_1+k_2}),
\end{equation}
o\`u $*$ d\'esigne comme ci-dessus l'\'etoile-produit associatif
$m+\Cal U(\hbar \gamma)$. La compatibilit\'e de ce produit avec
l'op\'erateur de cobord $[-,*]$ est imm\'ediate, et montre que le
produit $\cup$ induit un produit associatif $\cup$ sur sur l'espace
de cohomologie $H_{\hbar\tilde\gamma}$ du deuxi\`eme espace
tangent\footnote {Par rapport \`a \cite {K} chap. 8, nous nous pla\c
cons dans le cas simplifi\'e o\`u $\mg m$ est suppos\'ee
concentr\'ee en degr\'e z\'ero.}.
\begin{theoreme}\label{principal}
Dans le cas particulier $M=\R^d$ et
o\`u $\Cal U$ est le $L_\infty$-morphisme donn\'e dans \cite {K} \S\
6.4, la d\'eriv\'ee $\Cal U'_{\hbar\gamma}$ induit un isomorphisme
d'alg\`ebres de l'espace de cohomologie $H_{\hbar\gamma}$ de
l'espace tangent $T_{\hbar\gamma}\g g_1[1]$ sur l'espace de
cohomologie $H_{\hbar\tilde\gamma}$ de l'espace tangent
$T_{\hbar\tilde\gamma}\g g_2[1]$. Autrement dit pour tout couple
$(\alpha,\beta)$ de multi-champs de vecteurs tels que
$[\alpha,\gamma]=[\beta,\gamma]=0$, on a~:
\begin{equation}\label{multiplicativite}
\Cal U'_{\hbar\gamma}(\alpha\cup\beta)=\Cal U'_{\hbar\gamma}(\alpha)\cup\Cal
 U'_{\hbar\gamma}(\beta)+D,
\end{equation}
o\`u $D$ est un cobord de Hochschild pour l'alg\`ebre d\'eform\'ee
$(C^\infty(M)[[\hbar]],*)$.
\end{theoreme}
Cet article r\'edig\'e d\'ebut 2001 (arXiv : math/QA/0106205) est
consacr\'e \`a la d\'emons\-tration d\'etaill\'ee de ce th\'eor\`eme
de Kontsevich. Elle repose sur un argument d'homotopie qui est
esquiss\'e dans \cite K au paragraphe 8.1, et d\'evelopp\'e dans
\cite {ADS} dans le cas particulier o\`u $\alpha$ et $\beta$ sont de
degr\'e minimal $-2$ dans $\g g_1[1]$, c'est-\`a-dire appartiennent
\`a $C^\infty(M)[[\hbar]]$. Le terme de cobord $D$ n'appara\^\i t
pas dans ce cas. Il nous paraissait important de pr\'eciser les
\'el\'ements de cette d\'emonstration vu l'importance du r\'esultat
et la complexit\'e des m\'ethodes imagin\'ees par Kontsevich. Nous
mettons notamment en valeur une composante de bord nouvelle par
rapport \`a l'article \cite {ADS}, nous pr\'ecisions les
orientations des vari\'et\'es de bord qui produisent les termes de
co-bord et clarifions certaines conventions combinatoires source de
confusion. Cet article se veut \'el\'ementaire et p\'edagogique dans
son approche. Signalons aussi les travaux de Mochizuki \cite {Mo}
sur les structures $\Cal A^\infty$ associ\'ees aux structures
tangentes et de Tamarkin \cite {Ta} sur les structures $\Cal
G^\infty$. Ces articles moins \'el\'ementaires expliquent pourquoi
il est "naturel" d'avoir un transport de la structure associative en
cohomologie (voir aussi \cite {GH} sur les structures de Tamarkin
tangentes). Pour un analogue partiel du th\'eor\`eme 1.2 dans le cas
d'une vari\'et\'e quelconque on se reportera \`a \cite {CFT}.
\paragraph{Conventions sur le produit ext\'erieur :}
Soit $V$ un espace vectoriel de dimension finie. On identifie
$x_1\wwedge x_k$ avec la projection~:
$${1\over k!}\sum\limits_{\sigma\in S_k}\varepsilon(\sigma)x_{\sigma_1}\otimes\cdots\otimes
x_{\sigma_k}$$ de $x_1\otimes\cdots\otimes x_k$ sur les tenseurs
antisym\'etriques. On d\'efinit le couplage entre les puissances
ext\'erieures de $V$ et de son dual $V^*$ par la formule~:
\begin{equation}\label{ext}
<x_1\wwedge x_k,\, \xi_1\wwedge \xi_k> ={1\over
k!}\sum\limits_{\sigma\in S_k}\varepsilon(\sigma)<x_{\sigma
1},\,\xi_1>\cdots <x_{\sigma k},\,\xi_k>.
\end{equation}
L'article
est organis\'e de la fa\c con suivante. La partie 2 est
consacr\'ee \`a la description de la formule de Kontsevich. La
partie 3 traite de mani\`ere pr\'ecise et compl\`ete  l'argument
d'homotopie. La partie 4 rassemble les \'el\'ements pour la
d\'emonstration du th\'eor\`eme 1.2.
\section{Rappels sur la construction du quasi-isomorphisme L$_\infty$}
\setcounter{equation}{0}

 Nous rassemblons ici les faits relatifs
aux graphes et aux poids introduits par M. Kontsevich dans \cite
K, qui permettent de construire explicitement le
$L_\infty$-quasi-isomorphisme $\Cal U$ dans le cas o\`u la
vari\'et\'e est $\R^d$, et de montrer qu'il poss\`ede les
propri\'et\'es attendues. On se reportera \`a \cite K chap. 5 et 6
pour les d\'etails, ainsi qu'\`a \cite {AMM} pour la d\'elicate
question du choix des signes.
\subsection{Graphes et poids}
On se place dans le cas plat $M=\R^d$. Le
$L_\infty$-quasi-isomorphisme $\Cal U$ est uniquement d\'etermin\'e
par ses coefficients de Taylor~:
$$\Cal U_n:S^n(\g g_1[1])\longrightarrow \g g_2[1].$$
Si les $\alpha_k$ sont des $s_k$-champs de vecteurs, ils sont de
degr\'e $s_k-2$ dans l'espace d\'ecal\'e $\g g_1[1]$, et donc $\Cal
U_n(\alpha_1\cdots \alpha_n)$ est d'ordre $s_1+\cdots +s_n-2n$ dans
$\g g_2[1]$. C'est donc un op\'erateur $m$-diff\'erentiel, avec~:
\begin{equation}\label{magique}
\sum\limits_{k=1}^n s_k=2n+m-2.
\end{equation}
Les coefficients de Taylor sont construits \`a l'aide de poids et de
graphes~: on d\'esigne par $G_{n,m}$ l'ensemble des graphes
\'etiquet\'es et orient\'es ayant $n$ sommets du premier type
(sommets a\'eriens) et $m$ sommets du deuxi\`eme type (sommets
terrestres) tels que~:

\smallskip
1). Les ar\^etes partent toutes des sommets a\'eriens.\smallskip

2). Le but d'une ar\^ete est diff\'erent de sa source (il n'y a
pas de boucles).
\smallskip

3). Il n'y a pas d'ar\^etes multiples.

\ssq Par graphe \'etiquet\'e on entend un graphe $\Gamma$ muni
d'un ordre total sur l'ensemble $E_\Gamma$ de ses ar\^etes,
compatible avec l'ordre des sommets. A tout graphe \'etiquet\'e
$\Gamma\in G_{n,m}$, et \`a tout $n$-uple de multi-champs de
vecteurs $\uple \alpha n$ on peut associer de mani\`ere naturelle
un op\'erateur $m$-diff\'erentiel $B_\Gamma
(\alpha_1\otimes\cdots\otimes \alpha_n)$ lorsque pour tout
$j\in\{1,\ldots ,n\}$, $\alpha_j$ est un $s_j$-champ de vecteurs,
o\`u $s_j$ d\'esigne le nombre d'ar\^etes qui partent du sommet
a\'erien num\'ero $j$ (\cite {K} \S~ 6.3). \ssq L'op\'erateur
$B_\Gamma (\alpha_1\otimes\cdots\otimes \alpha_n)$ est construit
de la fa\c con suivante~: soit $E_\Gamma$ l'ensemble des ar\^etes
de $\Gamma$. On d\'esigne par $\{e_k^1,\ldots,e_k^{s_k}\}$ le
sous-ensemble ordonn\'e de $E_\Gamma$ des ar\^etes partant du
sommet a\'erien $k$. A toute application $I:E_\Gamma\to
\{1,\ldots, d\}$ et \`a tout sommet du graphe $\Gamma$ (de type
a\'erien ou terrestre) on associe l'op\'erateur diff\'erentiel \`a
coefficient constant~:
\begin{equation}
D_{I(x)}=\prod_{e=(-,x)}\bord_{I(e)},
\end{equation}
o\`u pour tout $i\in\{1,\ldots,d\}$ on d\'esigne par $\bord_i$
l'op\'erateur de d\'erivation partielle par rapport \`a la
$i^{\hbox{\sevenrm \`eme}}$ variable. Le produit est pris pour
toutes les ar\^etes qui arrivent au sommet $x$. Soit pour tout
sommet a\'erien $k$, pour tout $s_k$-champ de vecteurs $\gamma_k$
(o\`u $s_k$ est le nombre d'ar\^etes partant de $k$) et pour toute
application $I:E_\Gamma\to \{1,\ldots, d\}$ le coefficient~:
\begin{eqnarray}
\gamma_k^I   =\gamma_k^{I(e_k^1)\cdots I(e_k^{s_k})} &=
<\gamma_k,\,dx_{I(e_k^1)}\wwedge dx_{I(e_k^{s_k})}>\\\nonumber
        &=<\gamma_k,\,dx_{I(e_k^1)}\otimes\cdots\otimes dx_{I(e_k^{s_k})}>.
\end{eqnarray}
On pose alors~:
\begin{equation}\label{multidif}
\Cal B_\Gamma(\gamma_1\otimes\cdots\otimes
\gamma_n)(f_1\otimes\cdots\otimes f_n)= \sum\limits_{I:E_\Gamma\to
\{1,\ldots, d\}}\prod_{k=1}^n D_{I(k)}\gamma_k^I
        \prod_{l=1}^mD_{I(\overline l)}f_l.
\end{equation}
\smallskip
Le coefficient de Taylor $\Cal U_n$ est alors donn\'e par la
formule~:
\begin{equation}\label{taylor}
\Cal U_n(\alpha_1\cdots \alpha_n)=\sum\limits_{\Gamma\in G_{n,m}}
        W_\Gamma \Cal B_\Gamma(\alpha_1\otimes\cdots\otimes\alpha_n),
\end{equation}
o\`u l'entier $m$ est reli\'e \`a $n$ et aux $\alpha_j$ par la
formule \ref{magique} ci-dessus. On \'ecrit aussi~:
\begin{equation}\label{taylor2}
\Cal U_n=\sum\limits_{G_n}W_\Gamma \Cal B_\Gamma,
\end{equation}
o\`u $G_n$ d\'esigne l'ensemble des graphes admissibles \`a $n$
sommets a\'eriens.
\ssq
Le poids $W_\Gamma$ est nul sauf si le
nombre d'ar\^etes $|E_\Gamma|$ du graphe $\Gamma$ est
pr\'ecis\'ement \'egal \`a $2n+m-2$. Il s'obtient en int\'egrant une
forme ferm\'ee $\omega_\Gamma$ de degr\'e $|E_\Gamma|$ sur une
composante connexe de la compactification de Fulton-McPherson d'un
espace de configuration qui est pr\'ecis\'ement de dimension
$2n+m-2$ (\cite {FM}, \cite {K} \S~5). Il d\'epend lui aussi d'un
ordre sur l'ensemble des ar\^etes, mais le produit
$W_\Gamma. \Cal B_\Gamma$ n'en d\'epend plus.
\ssq
On d\'esigne
par $\mop{Conf}_{n,m}$ l'ensemble des $(\uple pn,\uple qm)$ o\`u les
$p_j$ sont des points distincts appartenant au demi-plan de
Poincar\'e~:
$$\Cal H=\{z\in\C, \mop{Im}z>0\},$$
et o\`u les $q_j$ sont des points distincts sur $\R$ vu comme le
bord de $\Cal H$. Le groupe~:
$$G=\{z\mapsto az+b \hbox{ avec }(a,b)\in\R \hbox{ et }a>0\}$$
agit librement sur $\mop{Conf}_{n,m}$. Le quotient~:
$$C_{n,m}=\mop{Conf}_{n,m}/G$$
est une vari\'et\'e de dimension $2n+m-2$. On a une action naturelle
des groupes de permutations $S_n$ et $S_m$ sur $C_{n,m}$, de sorte
qu'on peut parler de $C_{A,B}$ o\`u $A$ et $B$ sont deux ensembles
finis. La compactification de Fulton-MacPherson $\overline
{C_{A,B}}$ de $C_{A,B}$ est une vari\'et\'e \`a coins qui peut se
voir comme l'adh\'erence du plongement de $C_{A,B}$ dans une
vari\'et\'e compacte de grande dimension, produit de tores
$T=\R/2\pi\Z$ et d'espaces projectifs $\C P^2$ (\cite {AMM} \S\ I.2)~:
$$(\uple r{n+m})\longmapsto \bigl(\mop{Arg}(r_i-r_j)_{(i,j)},\,
        [r_i-r_j:r_j-r_k:r_k-r_i]_{(i,j,k)}\bigr).$$
Pour tout graphe $\Gamma\in G_{n,m}$ on construit une fonction
d'angle~:
$$\Phi_\Gamma:\overline C_{n,m}\longrightarrow T^{|E_\Gamma|}$$
de la fa\c con suivante~: on trace le graphe dans $\overline \Cal H$
en reliant les sommets par des g\'eod\'esiques pour la m\'etrique
hyperbolique, et \`a chaque ar\^ete $e=(a,b)$ on associe l'angle
$\varphi_e$ que fait la demi-droite verticale issue de $a$ avec
l'ar\^ete $e$~:
\dessin{40mm}{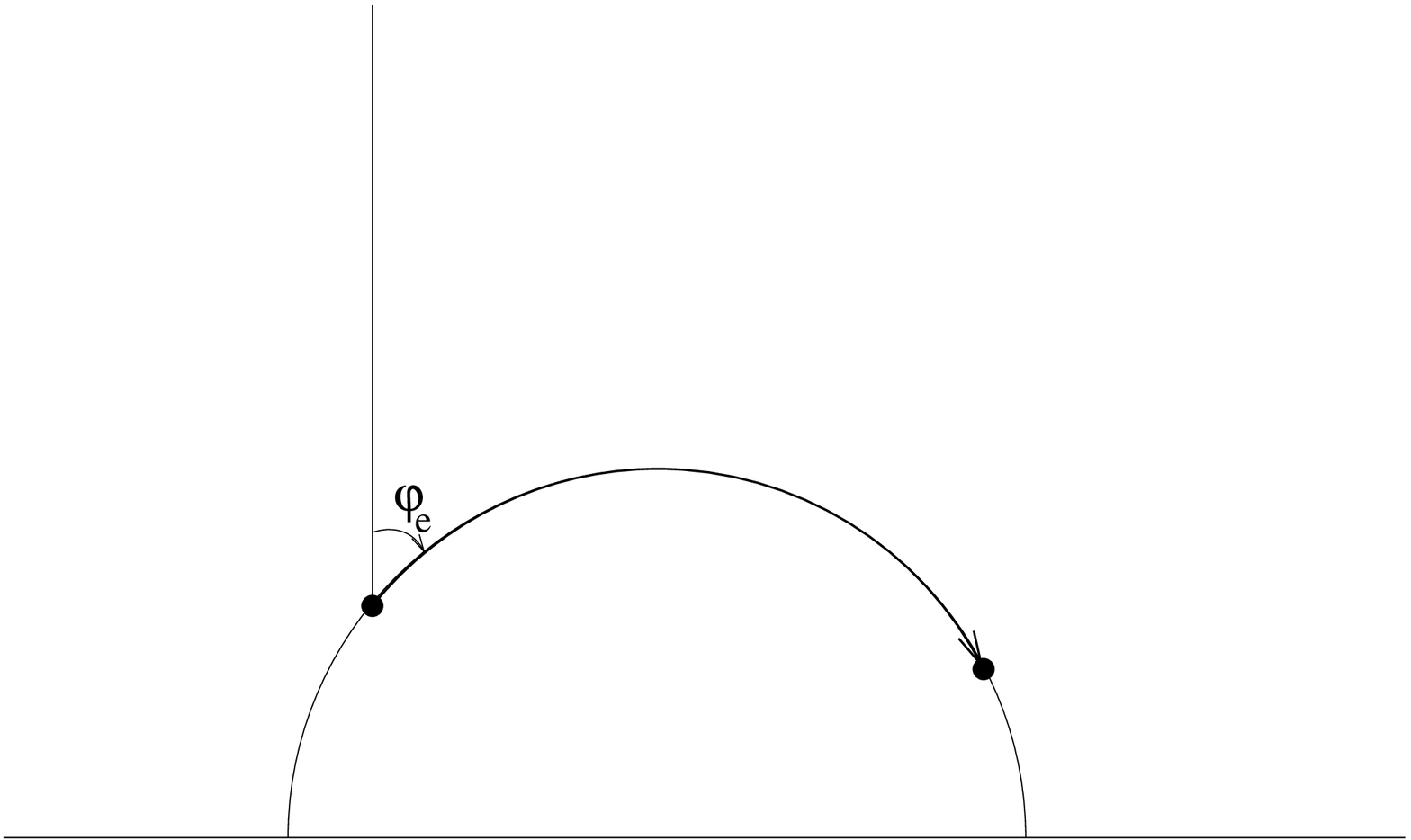}
En choisissant un ordre sur les ar\^etes ceci d\'efinit
$\Phi_\Gamma$ sur $C_{n,m}$ et on v\'erifie que cette application se
prolonge \`a la compactification. Soit $\omega_\Gamma$ la forme
diff\'erentielle $\Phi_\Gamma^*(dv)$ sur $C_{n,m}$ o\`u $dv$ est la
forme volume normalis\'ee sur $T^{|E_\Gamma|}$. La forme
$\omega_\Gamma$ se prolonge \`a la compactification. Soit
$C^+_{n,m}$ la composante connexe de $C_{n,m}$ o\`u les $\uple qm$
sont rang\'es par ordre croissant. Les orientations naturelles du
demi-plan $\Cal H$ et de $\R$ d\'efinissent une orientation de
$\mop{Conf}_{n,m}^+$, et par passage au quotient une orientation
naturelle de $C_{n,m}^+$, car l'action du groupe $G$ pr\'eserve
l'orientation. On d\'efinit alors le poids $W_\Gamma$ par~:
\begin{equation}\label{poids}
W_\Gamma=\int_{\overline C_{n,m}^+}\omega_\Gamma.
\end{equation}
{\sl Remarque\/}~: ce poids est un peu diff\'erent du poids
d\'efini par M. Kontsevich dans \cite {K} \S\ 6.2~: nous ne
multiplions pas l'int\'egrale par le facteur
$\bigl(\prod_{k=1}^n{1\over s_k!}\bigr)$. Ce facteur multiplicatif
est compens\'e par la convention sur les produits ext\'erieurs
(\ref{ext}) que nous avons adopt\'ee, qui diff\`ere de celle de
\cite {K} \S\ 6.3.
\subsection{Permutation des ar\^etes}
Soit $\Gamma$ un graphe admissible dans $G_{n,m}$. Le groupe
$S_{s_1}\times\cdots\times S_{s_n}$, produit des groupes de
permutations des ar\^etes attach\'es \`a chaque sommet, agit
naturellement sur $\Gamma$ par permutation de l'\'etiquetage des
ar\^etes. Il est clair que l'on a~:
\begin{eqnarray}
\Cal B_{\sigma.\Gamma}&=\varepsilon(\sigma)\Cal B_\Gamma \\
W_{\sigma.\Gamma}&=\varepsilon(\sigma)W_\Gamma,
\end{eqnarray}
de sorte que le produit $W_\Gamma.\Cal B_\Gamma$ ne d\'epend pas de
l'\'etiquetage.
\subsection{Stratification du bord}
Une suite de configurations $(D_n)$ tend vers un point $D$ du bord
de $\overline {C_{n,m}}$ lorsque des points se confondent ou tendent
vers la droite r\'eelle lorsqu'on choisit pour chaque configuration
un repr\'esentant $R_n$ en ``position standard'', c'est-\`a-dire
(par exemple), quitte \`a faire agir le groupe $G$, un
repr\'esentant tel que son diam\`etre soit \'egal \`a $1$ et tel que
l'abscisse de son barycentre soit $0$. \ssq On appelle {\sl
concentration\/} une suite $(R'_n)$ d'ensembles de points telle que
$R'_n\subset R_n$, telle que le cardinal de $R'_n$ soit
ind\'ependant de $n$, telle que la distance entre deux points
quelconques de $R'_n$ tende vers z\'ero, et telle que si $R'_n$ ne
contient qu'un seul point, alors ce point appartient \`a $\Cal H$
mais tend vers l'axe r\'eel. On appelle {\sl nuage\/} une
concentration maximale pour l'inclusion. On remarque que, par
d\'efinition de la position standard, un nuage ne peut pas contenir
tous les points de la configuration. \ssq On peut alors choisir un
nuage $(R'_n)$ , l'observer ``au microscope'', c'est-\`a-dire faire
agir un \'el\'ement du groupe $G$ (diff\'erent) sur chaque ensemble
de points $R'_n$ pour l'amener en position standard pour tout $n$,
et recommencer l'op\'eration ci-dessus. Le processus s'arr\^ete
apr\`es un nombre fini d'\'etapes, la derni\`ere \'etape ne
fournissant plus de nuages. Autrement dit le processus s'arr\^ete
lorsque tous les points sont distingu\'es individuellement et se
trouvent soit sur l'axe r\'eel, soit loin de l'axe r\'eel. \ssq Le
processus it\'eratif ci-dessus peut se d\'ecrire \`a l'aide d'un
arbre enracin\'e. Les sommets autres que les feuilles et la racine
correspondent au nuages, les feuilles correspondent aux points qui
apparaissent \`a la fin du processus. Un sommet {\sl noir\/}
correspond \`a un nuage a\'erien (ou \`a un point a\'erien si c'est
une feuille), et un sommet {\sl blanc\/} correspond \`a un nuage (ou
\`a un point) terrestre. La racine est par convention blanche (on
consid\`ere l'ensemble de la configuration vue de loin comme un
nuage terrestre), chaque sommet blanc a des descendants blancs ou
noirs, mais tous les descendants d'un sommet noir sont noirs. Chaque
arbre correspond \`a une strate du bord, dont la codimension est
\'egale au nombre $c(T)$ de sommets interm\'ediaires (ni racine ni
feuilles), ou si l'on pr\'ef\`ere, au nombre d'observations au
microscope n\'ecessaires. Voici un exemple repr\'esent\'e \`a l'aide
de nuages, avec l'arbre correspondant~:
\dessin{50mm}{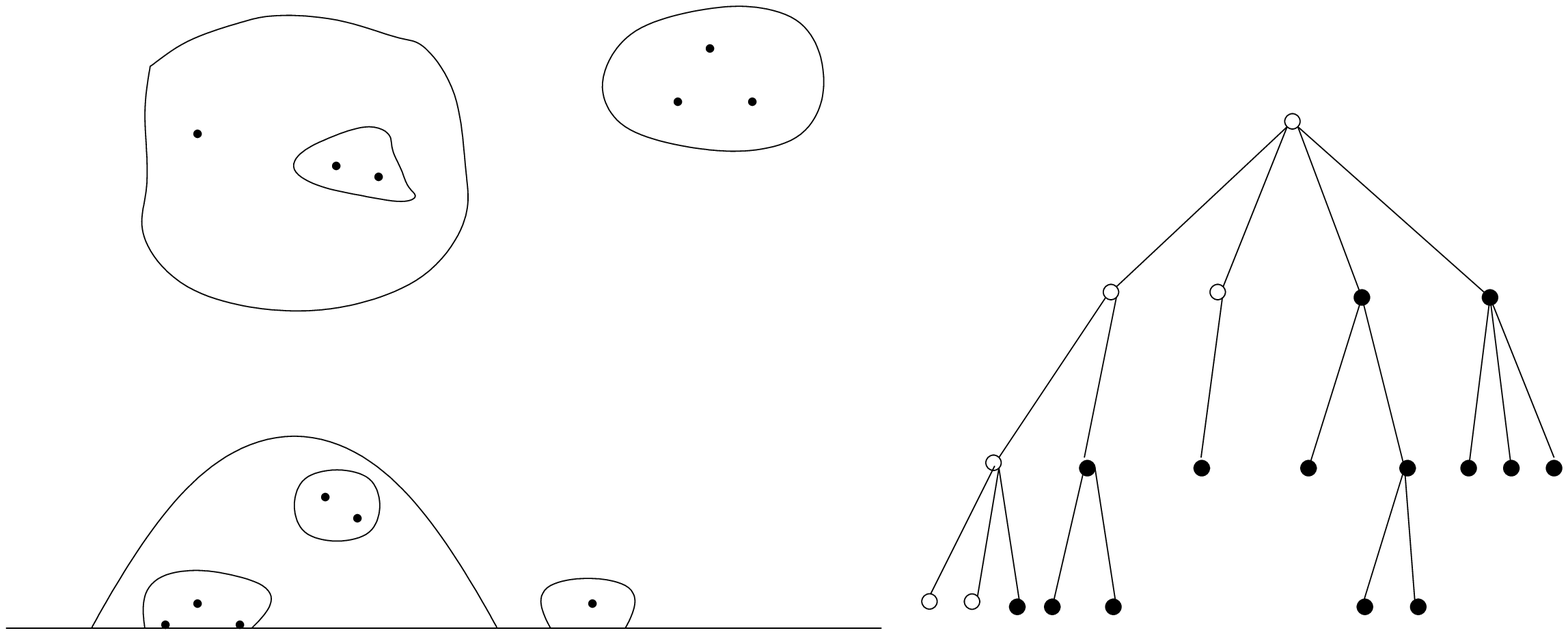}
Soit $\mop{Conf}_n$ l'ensemble des $n$-uplets de points distincts
dans $\C^n$. On note $C_n$ le quotient de $\mop{Conf}_n$ par
l'action du groupe $G'=\{z\mapsto az+b, \,a>0 \hbox{ et } b\in\C\}$.
C'est une vari\'et\'e de dimension $2n-3$, qu'il ne faut pas
confondre avec $C_{n,0}$. Le groupe de permutations $S_n$ agit
naturellement sur $C_n$, de sorte qu'on peut parler de $C_E$ o\`u
$E$ est un ensemble fini. Pour un arbre $T$ comme d\'ecrit ci-dessus
on peut d\'ecrire la strate correspondante $\bord_T C_{n,m}$ comme
un produit d'espaces de configuration $C_{A,B}$ et $C_E$ de la fa\c
con suivante~: on d\'esigne par $B$ l'ensemble des sommets blancs
qui ne sont ni feuille ni racine, et par $N$ l'ensemble des sommets
noirs qui ne sont ni feuille ni racine. Pour tout sommet $v$ qui
n'est pas une feuille (mais pouvant \^etre la racine) on d\'esigne
par $B_v$ l'ensemble de ses descendants blancs et par $N_v$
l'ensemble de ses descendants noirs. Alors on a~:
\begin{equation}
\bord_T C_{n,m}=\prod_{v \hbox{ \sevenrm blanc }}C_{N_v,B_v}
        \times \prod_{v \hbox{ \sevenrm noir }}C_{N_v}.
\end{equation}
on a $\displaystyle{|N|+n=\sum\limits_{v \hbox{ \sevenrm noir
}}|N_v|+\sum\limits_{v \hbox{ \sevenrm blanc }}|N_v|}$, et
$\displaystyle{|B|+m=\sum\limits_{v \hbox{ \sevenrm blanc
}}|B_v|}$. Sachant que $n$ est le nombre de feuilles noires et $m$
le nombre de feuilles blanches on a donc~:
\begin{eqnarray}\nonumber
\mop{dim}\bord_T     &=\sum\limits_{v \hbox{ \sevenrm blanc }}
(2|N_v|+|B_v|-2)+\sum\limits_{v \hbox{ \sevenrm noir
}}(2|N_v|-3)\\
        &=2n+2|N|+m+|B|-2(|B|+1)-3|N|\\\nonumber
        &=2n+m-2-|N|-|B|.
\end{eqnarray}
La codimension de $\bord_T$ dans $\overline{C_{n,m}}$ est donc
\'egale au nombre total $c(T)$ de sommets interm\'ediaires,
c'est-\`a-dire au nombre n\'ecessaire de grossissements au
microscope.
\subsection{Principe de d\'emonstration du th\'eor\`eme de
formalit\'e}
On reprend les notations de l'introduction : l'\'equation de
formalit\'e $\Cal U\circ Q^1=Q^2\circ \Cal U$ se d\'eveloppe en une
suite $(E_n)$ d'\'equations exprimant le coefficient de Taylor $\Cal
U_n$ en fonction de $Q_2^1$, $Q_1^2$, $Q_2^2$ et les $\Cal U_k$ pour
$k<n$. On montre, en rempla\c cant $\Cal U$ par son expression
donn\'ee au \S\ 2.1, que l'\'equation $(E_n)$ est \'equivalente
\`a~:
\begin{equation}
\sum\limits_{\Gamma\in G_n}C_\Gamma\Cal B_\Gamma=0
\end{equation}
o\`u pour $\Gamma\in G_{n,m}$ le poids $C_\Gamma$ est donn\'e par~:
\begin{equation}
C_\Gamma=\sum\limits_{T,c(T)=1}\int_{\bord_T \overline
{C_{n,m}}^+}\omega_{\Gamma}
        =\int_{\bord \overline{C_{n,m}}^+}\omega_\Gamma.
\end{equation}
On exprime ici que l'int\'egrale sur le bord est la somme des
int\'egrales sur les strates de codimension $1$ convenablement
orient\'ees (cf. \cite {AMM} chap. 1). Or $d\omega_\Gamma=0$. La
formule de Stokes sur la vari\'et\'e \`a coins
$\overline{C_{n,m}}^+$ implique l'annulation de tous les poids
$C_\Gamma$, et donc le th\'eor\`eme de formalit\'e (\cite {K} chap.
6, \cite{AMM}).
\section{L'argument d'homotopie}
\setcounter{equation}{0} On reprend les notations de
l'introduction. Il s'agit de comparer les deux quantit\'es $\Cal
U'_{\hbar\gamma}(\alpha\cup\beta)$ et $\Cal
U'_{\hbar\gamma}(\alpha)\cup \Cal U'_{\hbar\gamma}(\beta)$ en les
exprimant \`a l'aide de graphes et de poids. L'id\'ee est de
mettre en \'evidence une famille $Z_t$ de sous-vari\'et\'es \`a
coins de codimension 2 (pour $t\in[0,1]$) et une forme ferm\'ee
$\omega$ \`a valeurs dans les op\'erateurs multi-diff\'erentiels
telle que
$$\int_{Z_0}\omega = \Cal U'_{\hbar\gamma}(\alpha\cup\beta) \hbox to 20mm{\hfill et \hfill}
        \int_{Z_1}\omega        =\Cal U'_{\hbar\gamma}(\alpha)\cup \Cal U'_{\hbar\gamma}(\beta).$$
Gr\^ace \`a la formule de Stokes la diff\'erence entre les deux
quantit\'es s'exprime alors par l'int\'egrale sur $Y$ de $\omega$,
o\`u $Y$ d\'esigne la r\'eunion des bords des $Z_t$ pour $t\in
[0,1]$.
\subsection{L'oeil $\overline{C_{2,0}}$.}
Les diff\'erentes strates de $\overline{C_{2,0}}$ sont $C_{2,0}$,
$C_2$ (lorsque les deux points se rapprochent dans $\Cal H$), deux
copies de $C_{1,1}$ (lorsque l'un ou l'autre point se rapproche de
$\R$), et $C_{0,2}$ (lorsque les deux points se rapprochent de $\R$)
qui est constitu\'e par deux points. On trace un chemin $\xi(t)$
dans $\overline{C_{2,0}}$ tel que $\xi(0)\in C_2$ et $\xi(1)\in
C_{0,2}^+$~: autrement dit on part de deux points confondus dans
$\Cal H$ pour aboutir \`a deux points distincts sur $\R$.
\dessin{20mm}{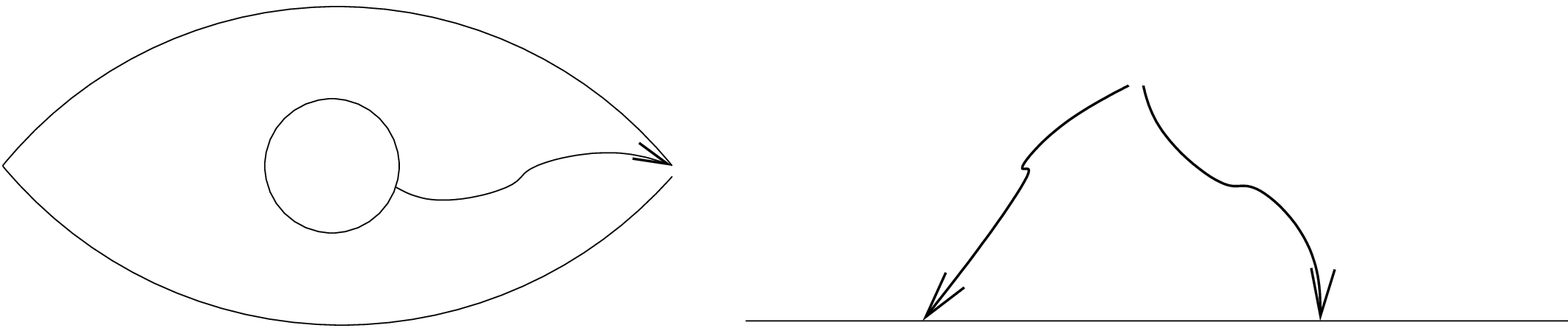}
\subsection{Une famille de sous-vari\'et\'es de codimension deux}
\ssq Soit $F:C_{n+2,m}^+\to C_{2,0}$ l'application surjective
consistant \`a oublier tous les points sauf les deux premiers
sommets a\'eriens. Cette application se prolonge naturellement par
continuit\'e en $F:\overline{C_{n+2,m}}^+\to \overline{C_{2,0}}$. On
d\'efinit alors pour $t\in[0,1]$~:
\begin{equation}\label{tranche}
Z_t=F\inver\bigl(\xi(t)\bigr).
\end{equation}
Il est clair que $Z_0$ et $Z_1$ sont contenus dans le bord de
$C_{n+2,m}^+$. Appelant $Z$ la r\'eunion des $Z_t$ pour $t\in]0,1[$,
et posant $Y=Z\cap \bord C_{n+2,m}^+$ on a
alors~\begin{equation}\label{boite} \bord Z=Z_0\amalg Z_1\amalg Y.
\end{equation}
Il faut maintenant d\'efinir pr\'ecis\'ement l'orientation des
espaces de configuration et des diff\'erentes strates qui
interviennent. Soit $\Omega$ la forme volume sur $C^+_{n+2,m}$
d\'efinie dans \cite {AMM} \S\ I.1, obtenue par passage au quotient
de la forme sur $\mop{Conf}_{n+2,m}$~:
$$\Omega_{\uple pn\, ; \,\uple qm}=dx_1\wedge dy_1\wwedge dx_{n+2}\wedge dy_{n+2}\wedge
        dq_1\wwedge dq_m$$
(avec $p_j=x_j+iy_j$).
\def\dn{d\mbox{\bf n}}
\dessin{50 mm}{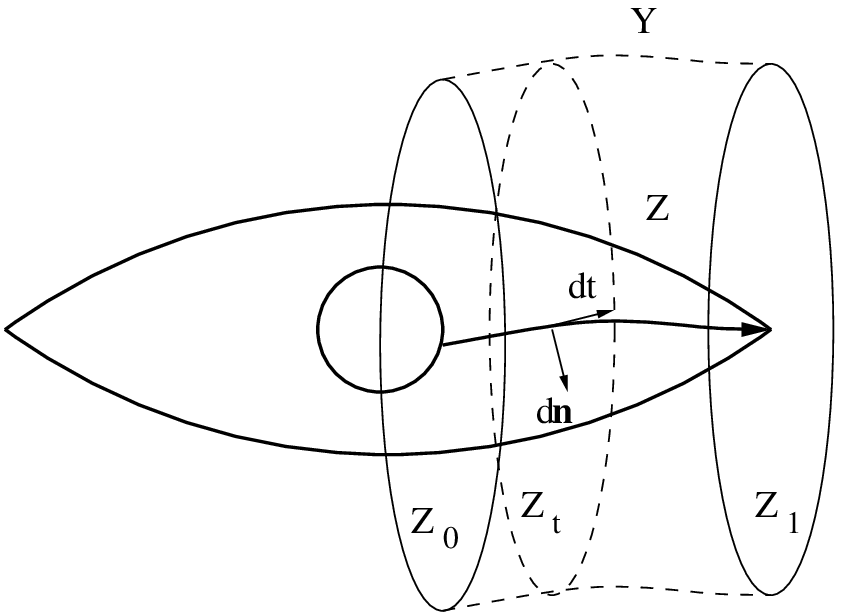}
Soit $\dn$ la normale au chemin $\xi$ dans l'oeil
$\overline{C_{2,0}}$ telle que $\dn\wedge dt$ forme un rep\`ere
direct du plan dans lequel est dessin\'e l'oeil (voir dessin). La
forme volume s'\'ecrit au voisinage d'un point quelconque de $Z_t$~:
$$\Omega=\dn\wedge dt\wedge \Omega',$$
ce qui d\'efinit une forme volume $\Omega'$ (et donc une
orientation) sur $Z_t$.
Soit $\Gamma$ un graphe admissible dans $G_{n+2,m}$, et
$\omega_\Gamma$ la forme diff\'erentielle associ\'ee. On pose~:
\begin{equation}
W_\Gamma^t=\int_{Z_t}\omega_\Gamma
\end{equation}
o\`u $Z_t$ est orient\'ee comme ci-dessus. Cette int\'egrale est
nulle sauf \'eventuellement si le nombre d'ar\^etes de $\Gamma$ est
\'egal \`a la dimension de $Z_t$, c'est-\`a-dire $2n+m$.
\begin{lemme}\label{poidsnul}
Le poids $W_\Gamma^0$ est nul si les sommets $1$ et $2$ sont
reli\'es par une ar\^ete, et si $1$ et $2$ ne sont pas reli\'es on a
avec nos choix d'orientation~:
\begin{equation}
W_\Gamma^0=W_{\Delta}=\int_{\overline{C_{n+1,m}}^+}\omega_{\Delta},
\end{equation}
o\`u $\Delta$ est le graphe de $G_{n+1,m}$ obtenu en confondant les
sommets $1$ et $2$.
\end{lemme}
\begin{preuve}
Il n'y a pas de strate de codimension $2$
contenue dans $Z_0$. On a donc~:
\begin{equation}
\mathop{Z_0}\limits^\circ=\coprod_{c(T)=1}Z_0\cap\bord_T C_{n+2,m}.
\end{equation}
D'apr\`es le \S\ I.3 les strates de codimension $1$ (les {\sl
faces\/}) sont obtenues en consid\'erant une seule fois un nuage de
points. Elles sont donc de deux types~: le type 1 pour un nuage
a\'erien, et le type 2 pour un nuage terrestre, qui contient donc
\`a la fois des points a\'eriens et des points terrestres. Il est
clair que seules les faces de type 1 ont une intersection non vide
avec $Z_0$, et que le nuage \`a consid\'erer contient alors
forc\'ement les sommets $1$ et $2$.
\ssq
Soit $A\subset\{1,\ldots
,n\}$ un tel nuage dans $C_{n+2,m}$, soit $\Sigma_A$ la face
associ\'ee, soit $\Gamma\subset G_{n+2,m}$ un graphe \`a $2n+m$
ar\^etes. Soit $\Gamma_A^2$ le graphe obtenu en contractant le nuage
$A$ sur un point, et soit $\Gamma_A^1$ le graphe interne, form\'e
des seuls points de $A$ et des ar\^etes qui relient deux points de
$A$. On a alors~:
\begin{equation}
\int_{\overline{\Sigma_A\cap Z_0}}\omega_\Gamma=\int_{\overline{C_A\cap Z_0}}
\omega_{\Gamma_A^1}
\int_{\overline{C_{n-|A|+3,m}}^+}\omega_{\Gamma_A^2}.
\end{equation}
Si les sommets $1$ et $2$ sont reli\'es par une ar\^ete $e$, l'angle
associ\'e $\varphi_e$ est constant sur $Z_0$, et donc le terme de
gauche du produit s'annule. Il suffit donc de consid\'erer les
graphes $\Gamma$ o\`u $1$ et $2$ ne sont pas reli\'es. Soit
$\wt\Gamma$ le graphe $\Gamma$ auquel on a rajout\'e (en premi\`ere
position pour l'ordre total sur les ar\^etes) une ar\^ete $e$
joignant $1$ \`a $2$. On a alors~:
\begin{equation}
\int_{\overline{C_A\cap
Z_0}}\omega_{\Gamma_A^1}=\int_{C_A}d\varphi_e\wedge\omega_{\Gamma_A^1}=
\int_{C_A}\omega_{\wt\Gamma_A^1}.
\end{equation}
Cette derni\`ere
int\'egrale est nulle d\`es que $|A|\ge 3$, d'apr\`es le lemme 6.6
de \cite K. On est donc ramen\'e au cas o\`u le nuage n'est
constitu\'e que des sommets $1$ et $2$. Le graphe interne est
r\'eduit \`a deux points reli\'es, et le graphe externe est le
graphe $\Delta$. On a alors~:
\begin{equation}
\int_{\overline{C_A\cap Z_0}}\omega_{\Gamma_A^1}
\int_{\overline{C_{n+1,m}}^+}\omega_{\Gamma_A^2}
        =\int_{C_2}d\varphi_e
\int_{\overline{C_{n+1,m}}^+}\omega_{\Delta} .
\end{equation}
Il n'y a donc qu'une seule strate $\Sigma_A$ de codimension 1 \`a
consid\'erer. On a donc, en utilisant l'orientation de la strate
donn\'ee par le lemme I.2.1 de \cite {AMM}~:
\begin{eqnarray}\nonumber
W_\Gamma^0=\int_{Z_0}\omega_\Gamma=\int_{\Sigma_A\cap
Z_0}\omega_\Gamma=-\int_{\Sigma_A}\omega_{\wt\Gamma}=\\-\big(-\int_{C_A}
d\varphi_e\int_{\overline{C_{n+1,m}^+}}\omega_{\Delta}\big)
=\int_{\overline{C_{n+1,m}^+}}\omega_{\Delta},
\end{eqnarray}
 ce qui d\'emontre le lemme. Le premier signe moins vient du fait que $Z_0$ est
 orient\'ee avec une normale rentrante, alors que les conventions de \cite {AMM}
  utilisent la normale sortante. il vient donc d'apr\`es \cite {AMM} lemme I.2.1
   un deuxi\`eme signe moins qui annule le premier.
\end{preuve}

Consid\'erons maintenant $Z_1$. Son intersection avec les strates de
codimension $1$ est vide, donc $Z_1$ est une r\'eunion de strates de
codimension $2$. Par d\'efinition de $Z_1$ ces strates sont des
$C_T$ o\`u $T$ est un arbre du type ci-dessous~:
\dessin{30mm}{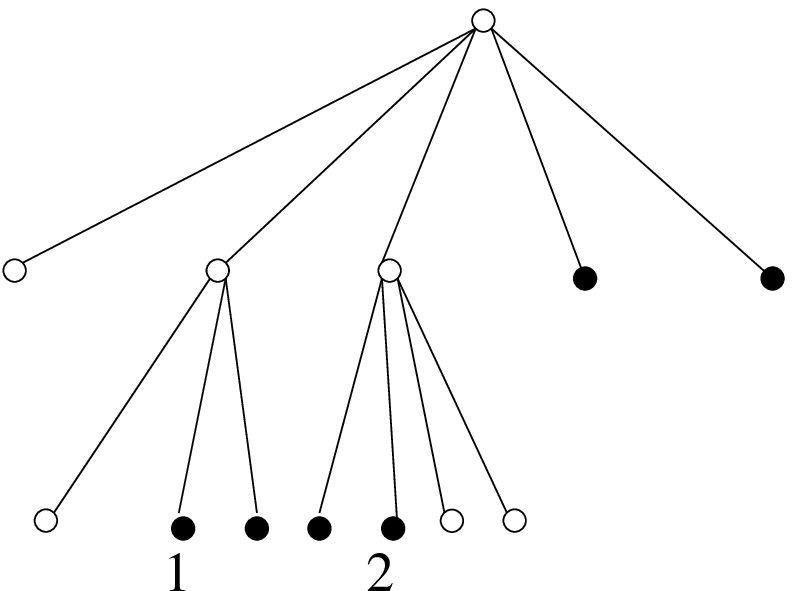}
c'est-\`a-dire qu'il y a deux nuages terrestres, un contenant le
sommet $1$ et l'autre contenant le sommet $2$. Une telle strate
s'\'ecrit donc~:
\begin{equation}\label{strate}
C_T=C_{n_1,m_1}\times C_{n_2,m_2}\times C_{n_3,m_3},
\end{equation}
avec $n_1+n_2+n_3=n+2$ et $m_1+m_2+(m_3-2)=m$. Son intersection avec
$Z_1$ est la composante connexe~:
$$C_T^+=C_{n_1,m_1}^+\times C_{n_2,m_2}^+\times C_{n_3,m_3}^+.$$
Les deux premiers facteurs repr\'esentent les deux nuages et le
dernier facteur repr\'esente la configuration ``avant
grossissement''. Soit $\Gamma$ un graphe dans $G_{n+2,m}$ contenant
$2n+m$ ar\^etes, soient $\Gamma_1^T$ et $\Gamma_2^T$ les deux
graphes internes correspondant aux deux nuages, et soit $\Gamma_3^T$
le graphe externe.
\begin{lemme}\label{trois}
On a $m_3=2$. Si $\Gamma$ poss\`ede une ou plusieurs fl\`eches
sortant du graphe interne $\Gamma_1$ ou $\Gamma_2$ alors
$W_\Gamma^1=0$. Dans le cas contraire on a~:
\begin{equation}
W_\Gamma^1=\sum\limits_T
W_{\Gamma_1^T}W_{\Gamma_2^T}W_{\Gamma_3^T},
\end{equation}
la somme portant sur tous les arbres $T$ d\'ecrits ci-dessus.
\end{lemme}
\begin{preuve}
Dans le premier cas, lorsqu'il existe une ``mauvaise fl\`eche''
$e$, l'angle $\varphi_e$ correspondant est constant sur $Z_1$, et
donc la forme diff\'erentielle \`a int\'egrer s'annule (\cite{K}
\S\ 6.4.2.2). Au signe pr\`es, la deuxi\`eme assertion d\'ecoule
naturellement de l'expression (\ref{strate}) donnant  $C_T$ comme
produit de trois espaces de configuration, et aussi du fait que
pour int\'egrer sur $Z_1$ il suffit de faire la somme des
int\'egrations sur les $C_T^+$ o\`u $T$ d\'ecrit l'ensemble des
arbres dessin\'es ci-dessus. Reste \`a pr\'eciser les orientations
pour d\'emontrer que le signe annonc\'e est le bon. \ssq
 On oriente chaque composante $C_{n_i,m_i}^+$ de la strate $C_T^+$ par sa forme
  volume $\Omega_i,i=1,2,3$. Le choix de normale que nous avons fait pour le chemin
  $\xi$ induit une orientation de chacune des paupi\`eres de l'oeil~: la paupi\`ere
  sup\'erieure est orient\'ee par la normale rentrante $\dn_1$, et la paupi\`ere
  inf\'erieure est ori\'ent\'ee par la normale sortante $\dn_2$. On s'en convainc
  en tra\c cant deux chemins $\xi_1$ et $\xi_2$ homotopes \`a $\xi$ dans
  $\overline{C_{2,0}}$, l'un se confondant avec la paupi\`ere sup\'erieure
  pour $t\ge 1/2$, et l'autre se confondant avec la paupi\`ere inf\'erieure pour $t\ge 1/2$~:
\dessin{30mm}{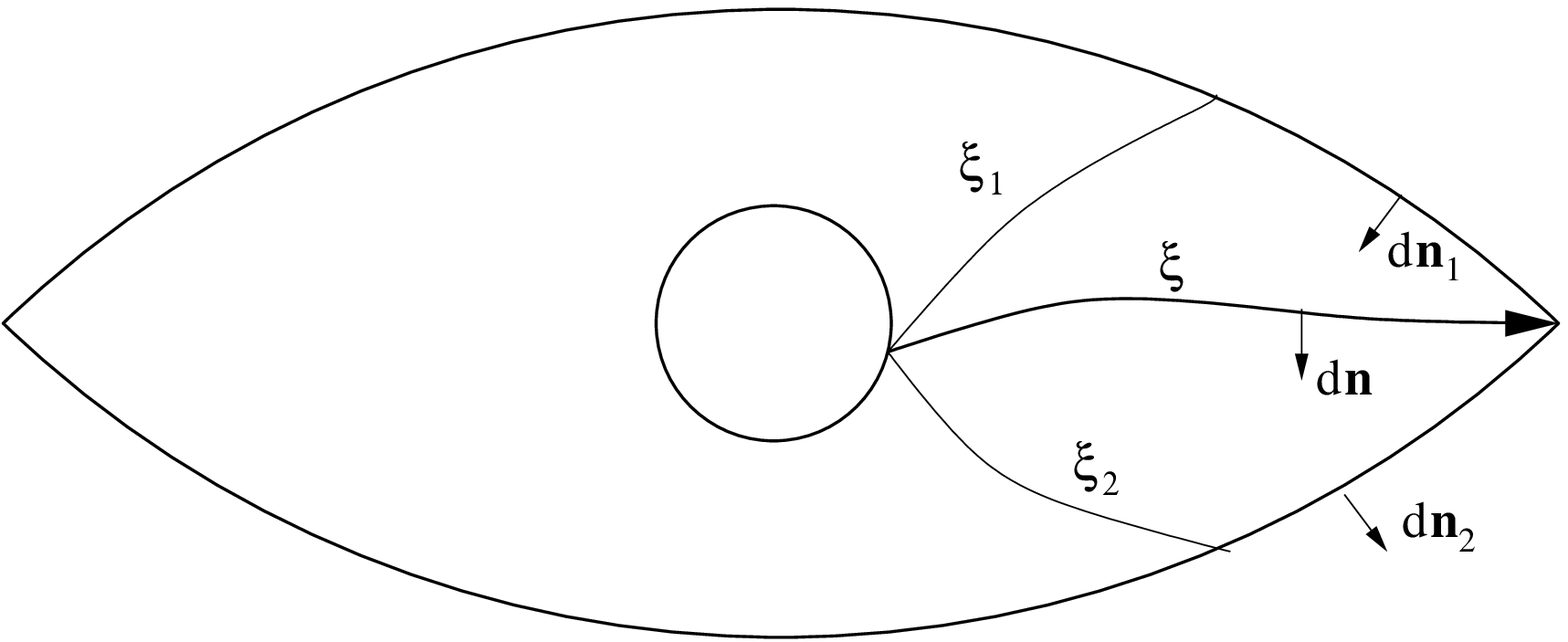}
\qquad De m\^eme que les coins de l'oeil forment l'intersection des
deux paupi\`eres, la strate $C_T$ de codimension 2 est
l'intersection des deux strates $C_{T_1}$ et $C_{T_2}$ de
codimension 1 obtenues respectivement en contractant seulement le
nuage de gauche (resp. de droite). D'apr\`es \cite {AMM} lemme I.2.2
la forme volume $\Omega$ s'\'ecrit~:
$$\Omega=-\dn_1\wedge\Omega(C_{T_1}^+)=(-1)^{lm_1+l+m_1+1}
\dn_1\wedge\Omega_1\wedge\Omega(C_{n_2+n_3,1+m_2}^+).$$
L'entier $l$
d\'esigne le nombre de points terrestres \`a gauche du nuage de
gauche~: on a donc $l=0$. La convention appliqu\'ee dans \cite {AMM}
est celle de la normale sortante~: comme $\dn_1$ est rentrante on a
un signe oppos\'e qui appara\^\i t par rapport \`a \cite {AMM} lemme
I.2.2. \ssq On contracte maintenant le nuage de droite. On a,
toujours gr\^ace \`a \cite {AMM} lemme I.2.2~:
$$\Omega(C_{n_2+n_3,1+m_2}^+)=\dn_2\wedge\Omega(C_{T'}^+)=
(-1)^{lm_2+l+m_2}\dn_2\wedge\Omega_2\wedge\Omega_3,$$ o\`u $C_{T'}$
d\'esigne la strate de codimension $1$ de $C_{n_2+n_3,1+m_2}$
hom\'eomorphe \`a $C_T$, obtenue en contractant le nuage de droite
une fois que le nuage de gauche est d\'ej\`a contract\'e. Cette
fois-ci $l=1$, car il y a un point \`a gauche du nuage~: celui
produit par la contraction du nuage de gauche. La normale $\dn_2$
\'etant sortante il n'y a pas de changement de signe ici. De plus on
peut manifestement remplacer $\dn_2$ par $dt$ pour peu que le chemin
$\xi$ arrive transversalement \`a la paupi\`ere inf\'erieure au coin
de l'oeil. On a donc finalement~:
\begin{eqnarray}
\Omega       &=(-1)^{m_1+1}(-1)^1\dn_1\wedge\Omega_1\wedge\dn_2
\wedge\Omega_2\wedge\Omega_3  \\\nonumber
                        &=\dn\wedge dt\wedge \Omega_1\wedge\Omega_2\wedge\Omega_3.
\end{eqnarray}
Chaque strate du type $T$ intervenant dans $Z_1$ est donc orient\'ee
par $\Omega_1\wedge\Omega_2\wedge\Omega_3$, d'o\`u le lemme.
\end{preuve}
\bsq

\begin{remarque}
Nous aurions pu choisir de contracter le nuage de
droite avant celui de gauche (autrement dit, nous aurions pu suivre
le chemin $\xi_2$ pour arriver \`a la strate $T$ au lieu de suivre
le chemin $\xi_1$). Nous laissons au lecteur le soin de v\'erifier
que cette alternative donne le m\^eme r\'esultat.
\end{remarque}
\subsection{Expression du cup-produit \`a l'aide de graphes et de
poids}\label{gp}
Les expressions suivantes constituent le premier pas de la
d\'emonstration du th\'eor\`eme 1.2
\begin{proposition}\label{importante}
Pour tout $\alpha,\beta\in\g g_1[1]$ on a~:
\begin{eqnarray}
&\Cal U'_{\hbar\gamma}(\alpha\cup\beta)
        =\sum\limits_{n\ge 0}{\hbar^n\over n!}\sum\limits\limits_{\Gamma\in
        G_{n+2,m}}
                W^0_\Gamma \Cal B_\Gamma(\alpha\otimes\beta\otimes\gamma
                        \otimes\cdots\otimes\gamma)\\
        &\Cal U'_{\hbar\gamma}(\alpha)\cup \Cal U'_{\hbar\gamma}(\beta)
        =\sum\limits_{n\ge 0}{\hbar^n\over n!}\sum\limits_{\Gamma\in G_{n+2,m}}
                W^1_\Gamma \Cal B_\Gamma(\alpha\otimes\beta\otimes\gamma
                        \otimes\cdots\otimes\gamma).
\end{eqnarray}
\end{proposition}

\begin{preuve}
On commence par d\'emontrer la premi\`ere assertion. On suppose
que $\alpha$ est homog\`ene de degr\'e $|\alpha|=m_1-1$ dans $\g
g_1$, et que $\beta$ est homog\`ene de degr\'e $|\beta|=m_2-1$
dans $\g g_2$. On rappelle que $\gamma$ est homog\`ene de degr\'e
$1$. L'expression explicite (\ref{derivee} de $\Cal
U'_{\hbar\gamma}$ donn\'ee dans l'introduction, ainsi que
l'expression explicite (\ref{taylor}) des $\Cal U_n$ \`a l'aide de
graphes et de poids, nous donnent~:
\begin{equation}
\Cal U'_{\hbar\gamma}(\alpha\cup\beta)=\Cal
U'_{\hbar\gamma}(\alpha\wedge\beta)= \sum\limits_{n\ge
0}{\hbar^n\over n!}\sum\limits_{\Gamma\in G_{n+1,m}}W_{\Gamma}\Cal
B_{\Gamma}
        \bigl((\alpha\wedge\beta)\otimes\gamma\otimes\cdots\otimes\gamma\bigr),
\end{equation}
avec $m=|\alpha|+|\beta|+2=m_1+m_2$. \ssq Introduisons
provisoirement, juste pour la suite de la d\'emonstration, une
classe sp\'e\-ciale de graphes~: on dira que l'ordre des ar\^etes
d'un graphe $\Gamma$ de $G_{n+1,m}$ est {\sl privil\'egi\'e par
rapport au premier sommet\/} si l'application qui \`a une ar\^ete
issue du premier sommet associe son sommet d'arriv\'ee est
croissante. Tout graphe $\Gamma\in G_{n+1,m}$ se d\'eduit d'un
graphe privil\'egi\'e par application d'une unique permutation
$\sigma$ des ar\^etes issues du premier sommet. On notera
$\varepsilon(\Gamma)$ la signature de cette permutation, et on
notera $\Gamma^{\hbox{\sevenrm priv}}$ le graphe privil\'egi\'e
naturellement associ\'e au graphe $\Gamma$. \ssq Soient $\Gamma$ et
$\Gamma'$ deux graphes admissibles (munis chacun d'un ordre sur
leurs ar\^etes compatible avec l'ordre sur leurs sommets). La
notation $\Gamma'\to\Gamma$ signifie que le graphe $\Gamma$ se
d\'eduit de $\Gamma'$ par compression d'un nuage de points. L'ordre
sur les ar\^etes de $\Gamma'$ impose l'ordre sur les ar\^etes de
$\Gamma$, et r\'eciproquement l'ordre sur les ar\^etes de $\Gamma$
impose l'ordre sur les ar\^etes de $\Gamma'$ qui sortent du nuage.
\ssq Les graphes $\Gamma'$ donnant $\Gamma$ par contraction peuvent
\^etre tr\`es diff\'erents des $\Gamma'$ donnant
$\Gamma^{\hbox{\sevenrm priv}}$ par contraction sur un sommet $s$~:
on s'en convainc (voir dessin ci-dessous) par un exemple simple dans
le cas o\`u $s$ n'a pas d'ar\^etes incidentes et o\`u le nuage est
compos\'e de deux sommets a\'eriens non reli\'es. Il n'y a qu'un
seul $\Gamma'$ par graphe contract\'e $\Gamma$ dans ce cas-l\`a.
\dessin{70mm}{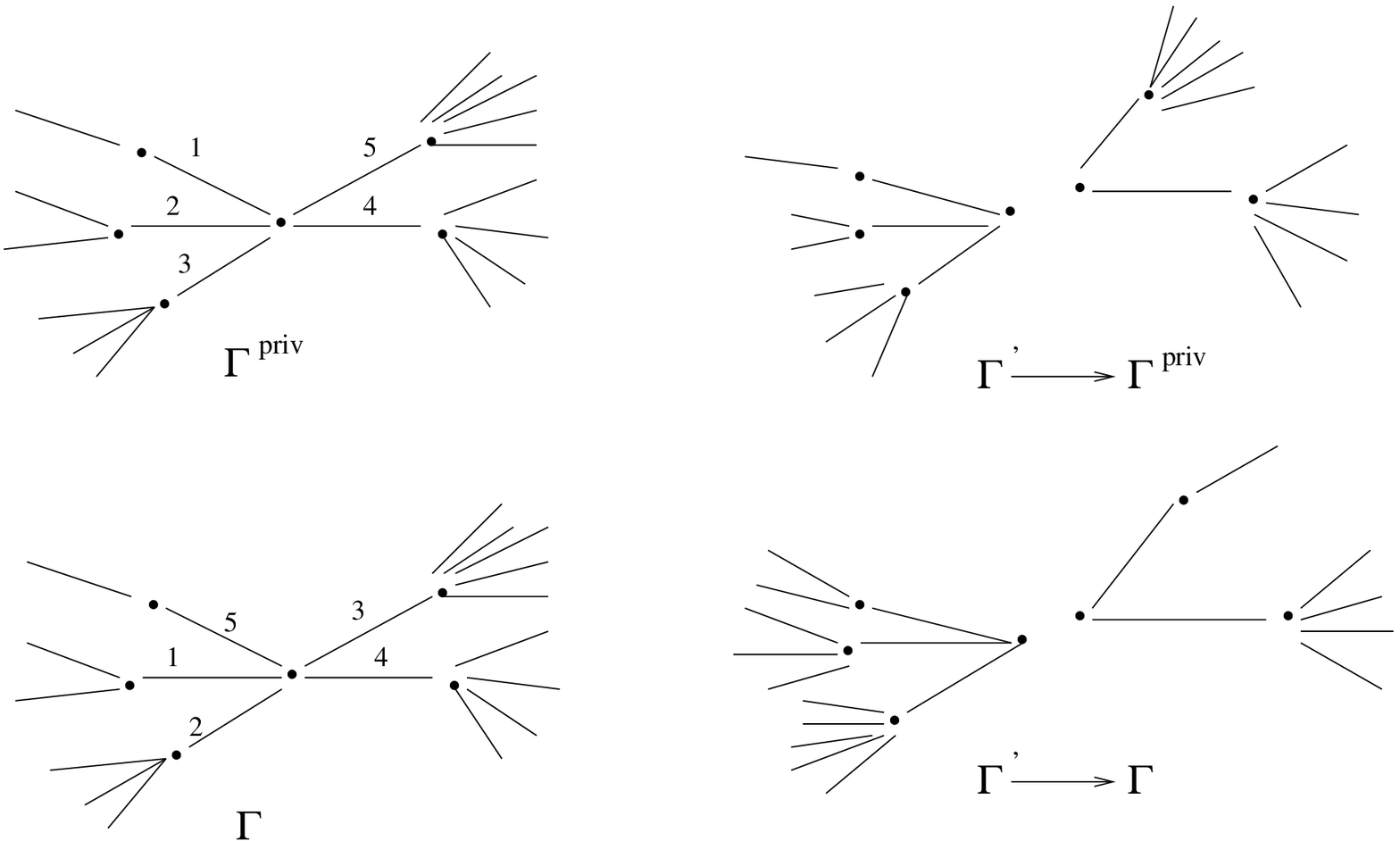}
\begin{lemme}
Soit $\Gamma_0$ un graphe dans $G_{n+1,m}$ privil\'egi\'e par
rapport au premier sommet. Alors,
\begin{eqnarray}\nonumber
(m_1+m_2)!\,\Cal B_{\Gamma_0}
        \bigl((\alpha& \hspace{-5cm}\wedge\beta)\otimes\gamma\otimes\cdots\otimes\gamma\bigr)
=\\\nonumber
=&\hspace{-2cm}\sum\limits_{\Gamma,\,\Gamma^{\hbox{\sevenrm
priv}}=\Gamma_0}
        \varepsilon(\Gamma)\Cal B_{\Gamma}
        \bigl((\alpha\wedge\beta)\otimes\gamma\otimes\cdots\otimes\gamma\bigr)\\
&=\sum\limits_{\Gamma,\,\Gamma^{\hbox{\sevenrm
priv}}=\Gamma_0}\varepsilon(\Gamma) \sum\limits_{\Gamma'\to\Gamma}
\Cal
B_{\Gamma'}(\alpha\otimes\beta\otimes\gamma\otimes\cdots\otimes\gamma),
\end{eqnarray}
o\`u la deuxi\`eme somme est prise sur tous les graphes $\Gamma'$
dans $G_{n+2,m}$ construits \`a partir de $\Gamma$ en d\'edoublant
le sommet $1$ de valence $m_1+m_2$ \'etiquet\'e par
$\alpha\wedge\beta$ en deux sommets $1$ et $2$ de valence $m_1$ et
$m_2$ respectivement, et en r\'epartissant les ar\^etes incidentes
entre les deux sommets de toutes les mani\`eres possibles.
\end{lemme}
\begin{preuve}
L'identit\'e bien connue sur les produits ext\'erieurs (avec la
convention (\ref{ext}))~: \begin{eqnarray*}
<\alpha\wedge\beta,\,dx_1\wwedge dx_{m_1+m_2}>={1\over
(m_1+m_2)!}\sum\limits_{\sigma\in S_{m_1+m_2}}
\end{eqnarray*}
\begin{eqnarray}\varepsilon(\sigma)
       <\alpha,\,dx_{\sigma_1}\wwedge dx_{\sigma_{m_1}}>
        <\beta,\,dx_{\sigma_{m_1+1}}\wwedge dx_{\sigma_{m_1+m_2}}>
\end{eqnarray}
 montre que l'on a pour toute application $I$ de
$E_{\Gamma'}$ dans $\{1,\ldots ,d\}$~:
\begin{equation}
(\alpha\wedge\beta)_1^I={1\over (m_1+m_2)!}\sum\limits_{\sigma\in
S_{m_1+m_2}}
        \varepsilon(\sigma)\alpha_1^{I\circ\sigma} \beta_2^{I\circ \sigma}.
\end{equation}
 On remplace alors $\Cal B_{\Gamma_0}
        \bigl((\alpha\wedge\beta)\otimes\gamma\otimes\cdots\otimes\gamma\bigr)$ par
         son d\'eveloppement explicite donn\'e par \ref{multidif}. L'utilisation de la formule
         de Leibniz fait alors appara\^\i tre une et une seule fois tous les graphes
          $\Gamma'\to\Gamma$ d\'ecrits ci-dessus (par r\'epartition des ar\^etes incidentes
          entre les sommets $1$ et $2$). Le groupe $S_{m_1+m_2}$ agit librement
          transitivement sur l'ensemble des graphes $\Gamma$ tels que
          $\Gamma^{\hbox{\sevenrm priv}}=\Gamma_0$, et le signe $\varepsilon(\Gamma)$
           est \'egal \`a la signature $\varepsilon(\sigma)$ o\`u $\sigma$ est
           la permutation telle que $\Gamma=\sigma\Gamma_0$.
\end{preuve}

Sachant que l'on a $W_\Gamma=\varepsilon(\Gamma)W_{\Gamma_0}$, on
d\'eduit du lemme pr\'ec\'edent et du lemme \ref{poidsnul}~:
\begin{equation}
\sum\limits_{\Gamma,\,\Gamma^{\hbox{\sevenrm priv}}=\Gamma_0}
        W_\Gamma\Cal B_{\Gamma}
        \bigl((\alpha\wedge\beta)\otimes\gamma\otimes\cdots\otimes\gamma\bigr)
=\sum\limits_{\Gamma,\,\Gamma^{\hbox{\sevenrm priv}}=\Gamma_0}
\sum\limits_{\Gamma'\to\Gamma}
        W_{\Gamma'}^0\Cal B_{\Gamma'}
        \bigl(\alpha\otimes\beta\otimes\gamma\otimes\cdots\otimes\gamma\bigr).
\end{equation}
On en d\'eduit la premi\`ere assertion de la proposition en sommant
sur tous les graphes $\Gamma$ de $G_{n+1,m}$ pour un $n$ donn\'e, en
multipliant par ${\hbar^n\over n!}$ et en sommant sur les entiers
$n$. \msq
\paragraph{Suite de la d\'emonstration de la proposition (deuxi\`eme assertion)}~:
 les \'egalit\'es (\ref{derivee}), (\ref{cup}) et (\ref{taylor}) nous donnent~:
\begin{eqnarray*}
&\Cal U'_{\hbar\gamma}(\alpha)\cup \Cal U'_{\hbar\gamma}(\beta)
        (f_1\otimes\cdots\otimes f_{m_1+m_2})=\sum\limits_{k,l\ge 0}{h^{k+l}\over k!l!}
        \sum\limits_{{\Gamma_1\in G_{k+1,m_1}\atop\Gamma_2\in G_{l+1,m_2}}}
W_{\Gamma_1}W_{\Gamma_2}\end{eqnarray*}
\begin{equation}
\Cal B_{\Gamma_1} (\alpha\otimes
\gamma\otimes\cdots\otimes\gamma)(f_1\otimes\cdots\otimes f_{m_1})
* \Cal B_{\Gamma_2} (\beta\otimes
\gamma\otimes\cdots\otimes\gamma)(f_{m_1+1}\otimes\cdots\otimes
f_{m_1+m_2}).
\end{equation}
D\'eveloppant le produit $*$ \`a l'aide de (\ref{tg2}) et
(\ref{taylor}) on obtient~:
\begin{eqnarray}\nonumber
&\Cal U'_{\hbar\gamma}(\alpha)\cup \Cal U'_{\hbar\gamma}(\beta)
        (f_1\otimes\cdots\otimes f_{m_1+m_2})=
\sum\limits_{k,l,r\ge 0}{h^{k+l+r}\over k!l!r!}
\sum\limits_{{{\Gamma_1\in G_{k+1,m_1}\atop\Gamma_2\in
G_{l+1,m_2}}\atop\Gamma_3\in G_{r,2}}}
W_{\Gamma_1}W_{\Gamma_2}W_{\Gamma_3}   \\\nonumber &\Cal
B_{\Gamma_3}(\gamma\otimes\cdots\otimes \gamma) \Bigl(\Cal
B_{\Gamma_1} (\alpha\otimes
\gamma\otimes\cdots\otimes\gamma)(f_1\otimes\cdots\otimes f_{m_1})
\,\otimes\\
&\otimes\,\Cal B_{\Gamma_2} (\beta\otimes
\gamma\otimes\cdots\otimes\gamma)(f_{m_1+1}\otimes\cdots\otimes
f_{m_1+m_2}) \Bigr).
\end{eqnarray}
Rappelons que toute strate $C_T$ de codimension $2$ incluse dans
$Z_1$ induit une d\'ecomposition d'un graphe $\Gamma\in G_{n+2,m}$
en deux graphes internes $\Gamma_1\in G_{n_1,m_1}$, $\Gamma_2\in
G_{n_2,m_2}$ et un graphe externe $\Gamma_3\in G_{n_3,m_3}$.
Compte tenu des degr\'es respectifs de $\alpha$, $\beta$ et
$\gamma$ les graphes $\Gamma$ intervenant dans le membre de droite
de la deuxi\`eme \'egalit\'e de la proposition \ref{importante}
v\'erifient $m_3=2$ pour toute d\'ecomposition de ce type :
autrement dit pour toute strate $C_T$ de codimension $2$ incluse
dans $Z_1$ les deux nuages absorbent tous les points terrestres.
R\'eciproquement, un triplet $(\Gamma_1,\Gamma_2,\Gamma_3)$ de
graphes comme dans la formule ci-dessus provient d'un ``grand''
graphe $\Gamma\in G_{n+2,m}$ avec $n=k+l+r$, $m=m_1+m_2$. \ssq
Reste le probl\`eme du d\'enombrement de tous les $\Gamma$ donnant
un triplet $(\Gamma_1,\Gamma_2,\Gamma_3)$ fix\'e~: on construit
$\Gamma$ \`a partir de ces trois composantes en choisissant la
mani\`ere de ``brancher'' $\Gamma_3$ sur $\Gamma_1$ et $\Gamma_2$,
et en choisissant un ordre sur les sommets a\'eriens de $\Gamma$
compatible avec l'ordre des sommets a\'eriens de $\Gamma_i$ pour
$i=1,2,3$. Il y a $\displaystyle |B_{k,l,r}|={(k+l+r)!\over
k!l!r!}$ choix possibles de tels ordres, et les choix de
branchements apparaissent tous par application de la r\`egle de
Leibniz lorsqu'on remplace les $\Cal B_{\Gamma_i}$ par leur
expression explicite donn\'ee par (\ref{multidif}). Le lemme
\ref{trois} permet alors de conclure. La proposition
\ref{importante} est donc d\'emontr\'ee.
\end{preuve}
\subsection{Le c\^ot\'e du cylindre : \'etude de $Y$.}\label{cylindre}
On peut voir le bord de $Z$ comme un cylindre born\'e (une
bo\^\i te de conserve) dont $Z_0$ serait le couvercle sup\'erieur,
$Z_1$ le couvercle inf\'erieur et $Y$ le c\^ot\'e. Il n'y a pas de
strates de codimension $2$ incluses dans $Y$, donc l'int\'erieur de
$Y$ est constitu\'e par la r\'eunion des $Z\cap \bord_T C_{n+2,m}$
avec $c(T)=1$, intersection de $Y$ avec les diff\'erentes strates de
codimension~$1$. On peut distinguer ici quatre types de strates,
suivant la position des sommets $1$ et $2$ par rapport au nuage~:
\dessin{100mm}{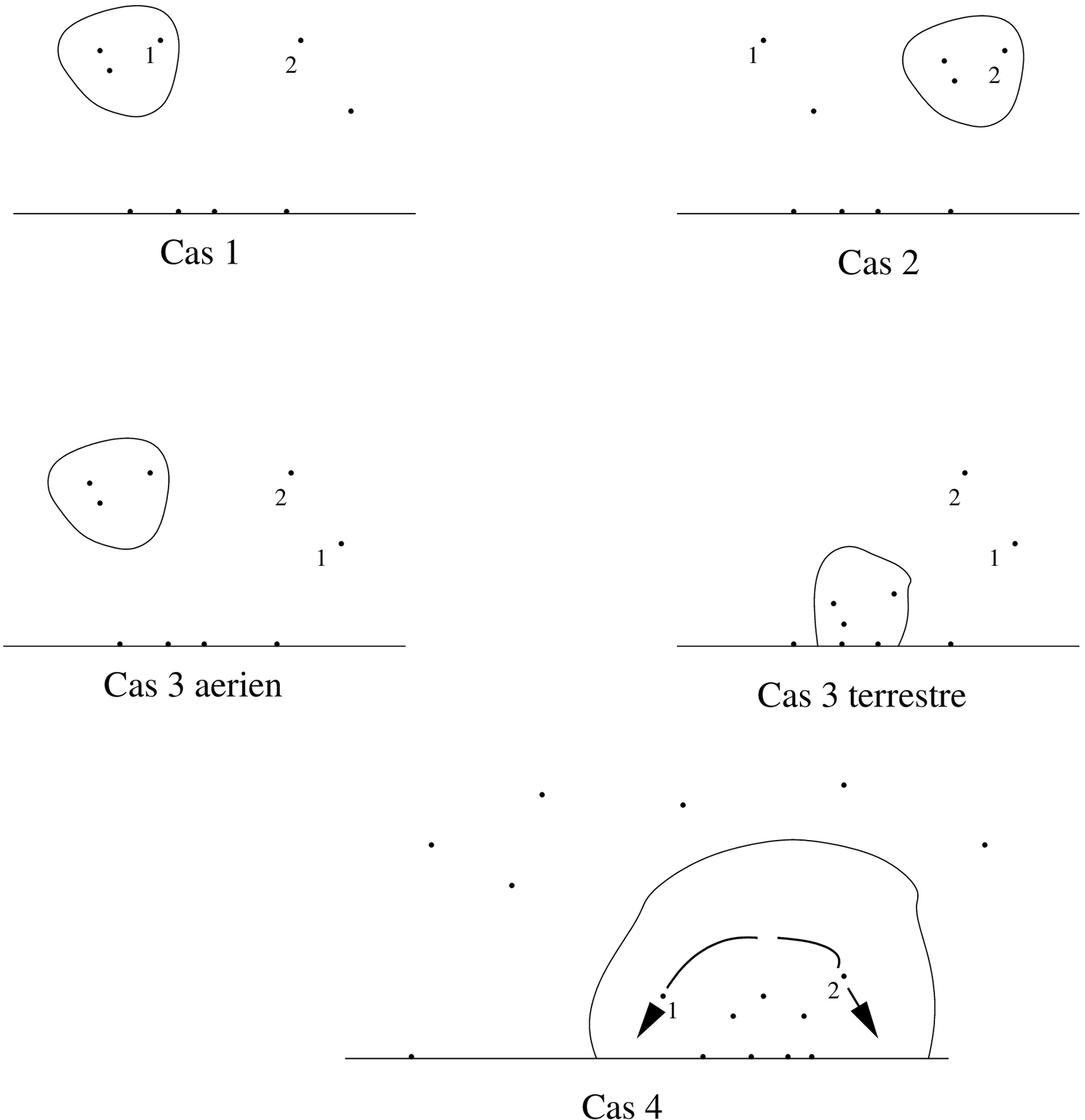}
{\sl Cas 1\/}~: le nuage contient le sommet $1$ et ne contient pas
le sommet $2$.

{\sl Cas 2\/}~: le nuage contient le sommet $2$ et ne contient pas
le sommet $1$.

{\sl Cas 3\/}~: le nuage ne contient ni le sommet $1$ ni le sommet
$2$.

{\sl Cas 4\/}~: le nuage contient \`a la fois les sommets $1$ et
$2$. \msq Dans les cas 1 et 2 le nuage est forc\'ement a\'erien, car
le chemin choisi dans l'oeil $C_{2,0}$ ne rencontre aucune des deux
copies de $C_{1,1}$. Dans le cas 3 le nuage peut \^etre a\'erien ou
terrestre. enfin dans le cas 4 le nuage est forc\'ement terrestre :
s'il \'etait a\'erien la strate serait contenue dans $F\inver (C_2)$
qui est d'intersection vide avec $Y$. Une telle strate contient
alors le chemin en entier, en ce sens que l'on a~:
$$C_T\cap Y=(F\restr{C_T})\inver\bigl(\xi([0,1])\bigr).$$
\break
\section{D\'emonstration du th\'eor\`eme 1.2.}
\setcounter{equation}{0}
\subsection{L'int\'egrale sur $Y$ comme cobord de Hochschild.}
Soit pour tout graphe $\Gamma$ dans $G_{n+2,m}$ le poids
$W''_\Gamma$ d\'efini par~:
\begin{equation}
W''_\Gamma=\int_Y \omega_\Gamma
\end{equation}
o\`u le bord $Y$ de $Z$ est orient\'e par la normale sortante.
D'apr\`es la formule de Stokes on a~:
\begin{equation}\label{stokes}
\int_{\bord Z}\omega_\Gamma=-W_\Gamma^0+W_\Gamma^1+W''_\Gamma=0.
\end{equation}
Le th\'eor\`eme d\'ecoulera donc imm\'ediatement de la proposition \ref{importante} et de la proposition suivante~:
\begin{proposition}\label{intermediaire}
Soit $\gamma$ un $2$-tenseur de Poisson formel, soit $*$
l'\'etoile-produit construit \`a partir de $\gamma$ \`a l'aide du
$L_\infty$-quasi-isomorphisme $\Cal U$, soient $\alpha$ un
$k_1$-champ de vecteurs et $\beta$ un $k_2$-champ de vecteurs.
Alors on a avec $m=k_1+k_2$~:
\begin{equation}\label{argh}
\sum\limits_{n\ge 0}{\hbar^n\over n!}\!\sum\limits_{\Gamma\in
G_{n+2,m}}
                W''_\Gamma \Cal B_\Gamma(\alpha\otimes\beta\otimes\gamma
                        \otimes\cdots\otimes\gamma)=\sum\limits_{n\ge 0}{\hbar^n\over n!}\!\sum\limits_{\Delta\in
G_{n+2,m-1}} \wt W_\Delta[*,\Cal
B_\Delta(\alpha\otimes\beta\otimes\gamma\otimes\cdots\otimes\gamma)]\end{equation}
$$
-\sum\limits_{n\ge 0}{\hbar^n\over n!}\sum\limits_{\Delta\in
G_{n+1,m}}\hskip -2mm \wt W_\Delta\Bigl(\Cal
B_\Delta([\alpha,\gamma]\otimes\beta\otimes
\gamma\otimes\cdots\otimes\gamma)+(-1)^{k_1}\Cal
B_\Delta(\alpha\otimes[\beta,\gamma]\otimes\gamma\otimes\cdots\otimes\gamma)\Bigr),
$$
o\`u $\wt W_\Delta$ d\'esigne l'int\'egrale de $\omega_\Delta$ sur
l'image r\'eciproque de $\xi([0,1])$ par l'application oubli
$F:\overline{C_{n+2,m-1}}\to \overline{C_{2,0}}$ ou
$F:\overline{C_{n+1,m}}\to \overline{C_{2,0}}$.
\end{proposition}
\begin{preuve}
On va montrer
ceci en d\'ecomposant $W''_\Gamma$ en somme de termes $W_\Gamma^T$
d'int\'e\-gra\-tions sur chacune des strates $C_T$ de codimension
$1$ rencontrant $Z$, qui sont de quatre types diff\'erents (voir \S\ \ref{gp}). On rappelle que, suivant \cite {AMM} \S IV.2, on doit pour
que le th\'eor\`eme de formalit\'e soit v\'erifi\'e remplacer le
crochet de Schouten par une petite modification de celui-ci,
d\'efinie par~:
\begin{equation}\label{schouten2}
[\gamma_1,\,\gamma_2]=-[\gamma_2,\,\gamma_1]_{\hbox{\sevenrm Schouten}}.
\end{equation}
C'est ce crochet modifi\'e qui est utilis\'e ici dans l'\'ecriture
des deux derniers termes de (\ref{argh}). Les paragraphes suivants
constituent la d\'emonstration de la proposition \ref{importante}.
On examine pr\'ecisemment la contribution de chacune des strates
d\'ecrites en \S \ref{cylindre}.
\subsection{L'int\'egration sur une strate de type 1 ou 2.}
\qquad On consid\`ere le type $1$, le type $2$ se traite de la
m\^eme mani\`ere. Le nuage est a\'erien et contient le sommet $1$
(\'etiquet\'e par $\alpha$). On peut supposer que le nuage ne
comporte que deux points, car pour un nombre de points sup\'erieur
l'int\'egrale de $\omega_\Gamma$ sur la strate correspondante
s'annule, d'apr\`es \cite {K} Lemma 6.6. On reprend la terminologie
et les notations du \S\ \ref{gp}.
\begin{lemme}\label{type1}
Soit $\Delta_0$ un graphe dans $G_{n+1,m}$ privil\'egi\'e par
rapport au premier sommet. Alors,

$$(k_1+k_2-1)!\Cal B_{\Delta_0}
        \bigl([\gamma_1,\,\gamma_2]\otimes\gamma_3
        \otimes\cdots\otimes\gamma_{n+2}\bigr)
        = \sum\limits_{\Delta,\,\Delta^{\hbox{\sevenrm priv}}=\Delta_0}
        \varepsilon(\Delta)\Cal B_{\Delta}
        \bigl([\gamma_1,\,\gamma_2]\otimes\gamma_3\otimes\cdots\otimes
        \gamma_{n+2}\bigr)$$\vspace{-0,5cm}
        \begin{equation}
=(-1)^{(k_1-1)k_2} \sum\limits_{\Delta,\,\Delta^{\hbox{\sevenrm
priv}}=\Delta_0}\varepsilon(\Delta) \sum\limits_{\Gamma\to\Delta}
\Cal B_{\Gamma}(\gamma_1\otimes\cdots\otimes\gamma_{n+2}),
\end{equation}
o\`u la deuxi\`eme somme est prise sur tous les graphes $\Gamma$
dans $G_{n+2,m}$ construits \`a partir de $\Delta$ en d\'edoublant
le sommet $1$ de valence $k_1+k_2-1$ \'etiquet\'e par
$[\gamma_1,\,\gamma_2]$ en deux sommets $1$ et $2$ de valence $k_1$
et $k_2$ respectivement, en tra\c cant une ar\^ete de $1$ vers $2$
ou de $2$ vers $1$, et en r\'epartissant les ar\^etes incidentes
entre les deux sommets de toutes les mani\`eres possibles.
\end{lemme}
\begin{preuve}
La formule g\'en\'erale pour le crochet de Schouten {\bf
modifi\'e} est la suivante~:
\begin{eqnarray*}
&[\xi_1\wwedge \xi_{k_1},\,\eta_1\wwedge \eta_{k_2}]
=-[\eta_1\wwedge \eta_{k_2},\,\xi_1\wwedge
\xi_{k_1}]_{\mbox{\sevenrm Schouten}}\\
&\hspace{-1cm}=\sum\limits_{r=1}^{k_1}\sum\limits_{s=1}^{k_2}(-1)^{r+s+1}[\eta_s,\xi_r]
        \wedge\eta_1\wwedge \wh{\eta_s}\wwedge \eta_{k_2}
        \wedge \xi_1
        \wwedge \wh{\xi_r}\wwedge \xi_{k_1}\end{eqnarray*}
        \begin{equation}
=\sum\limits_{r=1}^{k_1}\sum\limits_{s=1}^{k_2}(-1)^{r+s+(k_1-1)(k_2-1)}[\xi_r,\eta_s]\wedge
 \xi_1\wwedge \wh{\xi_r}\wwedge \xi_{k_1}
\wedge\eta_1\wwedge \wh{\eta_s}\wwedge \eta_{k_2}.
\end{equation}
On rappelle (\cite {AMM} \S\ II.4, IV.1 et IV.2) que le deuxi\`eme
coefficient de Taylor de la cod\'erivation $Q^1$ de $S^+(\g g_1[1])$
v\'erifie~:
\begin{equation}\label{antisym}
Q_2^1(\gamma_1.\gamma_2)=(-1)^{(k_1-1)k_2}[\gamma_1,\gamma_2],
\end{equation}
soit, avec $\gamma_1 =\xi_1\wwedge\xi_{k_1}$ et $\gamma_2
=\eta_1\wwedge \eta_{k_2}$~:

$$Q_2^1(\gamma_1.\gamma_2)=\sum\limits_{r=1}^{k_1}
\sum\limits_{s=1}^{k_2}(-1)^{r+s+k_1-1}[\xi_r,\eta_s]\wedge
 \xi_1\wwedge \wh{\xi_r}\wwedge \xi_{k_1}
\wedge\eta_1\wwedge \wh{\eta_s}\wwedge \eta_{k_2}.$$ On peut
supposer que $\xi_1=f\bord_{i_1}$, $\xi_r=\bord_{i_r}$ pour $r>1$,
$\eta_1=g\bord_{j_1}$ et $\eta_s=\bord_{j_s}$ pour $s>1$. On a
alors, toujours d'apr\`es \cite {AMM}\ \S\ IV.2~:
\begin{equation}\label{antisym2}
Q^1_2(\gamma_1.\gamma_2)=\gamma_1\bullet\gamma_2+(-1)^{k_1k_2}\gamma_2\bullet\gamma_1,
\end{equation}
avec~:
\begin{equation}\label{ptnoir}
\gamma_1\bullet\gamma_2
=\sum\limits_{r=1}^{k_1}(-1)^{r-1}f\bord_{i_r}g
\bord_{i_1}\wwedge\wh{\bord_{i_r}}\wwedge\bord_{i_{k_1}}\wedge\bord_{j_1}
\wwedge\bord_{j_{k_2}}.
\end{equation}
Le lemme d\'ecoulera donc
directement de la formule suivante~:
$$(k_1+k_2-1)!\Cal B_{\Delta_0}
        \bigl((\gamma_1\bullet\gamma_2)
        \otimes\gamma_3\otimes\cdots\otimes\gamma_{n+2}\bigr)=$$\vspace{-0,5cm}
        \begin{equation}\label{contract}
        \sum\limits_{\Delta,\,\Delta^{\hbox{\sevenrm priv}}=\Delta_0}
        \varepsilon(\Delta)\Cal B_{\Delta}
        \bigl((\gamma_1\bullet\gamma_2)\otimes\gamma_3\otimes\cdots\otimes\gamma_{n+2}\bigr)=\!\!\sum\limits_{\Delta,\,\Delta^{\mbox{\sevenrm
priv}}=\Delta_0}\varepsilon(\Delta) \sum\limits_{\Gamma\to\Delta}
\Cal B_{\Gamma}(\gamma_1\otimes\cdots\otimes\gamma_{n+2}),
\end{equation}
o\`u la somme interne porte sur les $\Gamma\to\Delta$ tels que la
fl\`eche reliant les sommets $1$ et $2$ est issue du sommet $1$.
On d\'eduit de (\ref{ptnoir}) la formule pour un coefficient
quelconque du multi-tenseur $\gamma_1\bullet\gamma_2$~:
\begin{eqnarray}\label{contract2}\nonumber
(\gamma_1\bullet\gamma_2)^{u_1\cdots u_{k_1+k_2-1}} ={1\over
(k_1+k_2-1)!}\sum\limits_{\sigma\in
S_{k_1+k_2-1}}\varepsilon(\sigma)
\sum\limits_{v=1}^d\sum\limits_{r=1}^{k_1}\\(-1)^{r-1}\gamma_1^{u_{\sigma_1}\cdots
u_{\sigma_{r-1}}vu_{\sigma_r}\cdots
u_{\sigma_{k_1-1}}}\bord_v(\gamma_2)^{u_{\sigma_{k_1}}\cdots
u_{\sigma_{k_1+k_2-1}}}.
\end{eqnarray}
Le groupe $S_{k_1+k_2-1}$ agit librement transitivement sur
l'ensemble des graphes $\Delta$ tels que $\Delta^{\hbox{\sevenrm
priv}}=\Delta_0$, et le signe $\varepsilon(\Delta)$ est \'egal \`a
la signature $\varepsilon(\sigma)$ o\`u $\sigma$ est la
permutation telle que $\Delta=\sigma\Delta_0$. La formule
(\ref{contract}) se d\'eduit alors de (\ref{antisym}),
(\ref{antisym2}) et (\ref{contract2}). Le lemme se d\'eduit
imm\'ediatement de (\ref{contract}) du fait que le signe
$(-1)^{k_1k_2}$ appara\^\i t lorsque l'on \'echange les deux
premiers sommets et qu'on r\'eordonne les ar\^etes de mani\`ere
compatible avec ce nouvel ordre sur les sommets.
\end{preuve}
\begin{corollaire}\label{cor}
Soit $C_T$ la strate de codimension $1$ de $\overline{C_{n+2,m}}^+$
correspondant au rapprochement des deux premiers sommets, orient\'ee
avec la convention de la normale sortante comme dans \cite {AMM} \S\
I.2.1. Soit $W_\Gamma^{T}$ l'int\'egrale de la forme $\omega_\Gamma$
sur la strate $C_T\cap Y$, o\`u $\Gamma$ est un graphe admissible.
On a alors~:
\begin{eqnarray}\nonumber
&\sum\limits_{\Delta,\,\Delta^{\hbox{\sevenrm priv}}=\Delta_0} \wt
W_\Delta\Cal B_\Delta
([\gamma_1,\gamma_2]\otimes\gamma_3\otimes\cdots\otimes
\gamma_{n+2})\\
&=(-1)^{(k_1-1)k_2}\sum\limits_{\Delta,\,\Delta^{\mbox{\sevenrm
priv}}=\Delta_0}
\sum\limits_{\Gamma\rightarrow\Delta}W_\Gamma^{T}\Cal
B_\Gamma(\gamma_1\otimes\cdots\otimes \gamma_{n+2}),
\end{eqnarray}
la somme \'etant prise sur tous les graphes
$\Gamma$ obtenus \`a partir de $\Delta$ en d\'edoublant le premier
sommet et en tra\c cant une ar\^ete d'un sommet vers l'autre.
\end{corollaire}
\begin{preuve}
La strate $C_T$ correspond au grossissement au microscope d'un
nuage a\'erien. On peut supposer que le nuage n'est compos\'e que de
deux points $a$ et $b$, en vertu du lemme 6.6 de \cite K. l'un de
ces sommets est le sommet $1$ ou $2$, l'autre est un sommet
diff\'erent de $1$ et $2$. On a~:
$$C_T\sim C_2\times C_{n+1,m},$$
et~:
$$C_T\cap Y\sim C_2\times \wt Z,$$
o\`u $\wt Z$ est l'image r\'eciproque de $\xi(]0,1[)$ par
l'application oubli de $C_{n+1,m}$ dans $C_{2,0}$. La strate $C_T$
est orient\'ee par $-\Omega_1\wedge\Omega_2$, d'apr\`es \cite {AMM}
\S\ I.2.1. On a donc~:
$$\int_{C_T\cap Y}\omega_{\Gamma}=-\int_{C_2}\omega_{\Gamma_1}\int_{\wt Z}\omega_{\Delta},$$
o\`u $\Gamma_1$ est le graphe interne et $\Delta$ le graphe externe.
Cette int\'egrale s'annule sauf \'eventuelle\-ment si les deux
points $a$ et $b$ sont reli\'es par une ar\^ete. Dans ce cas on a
donc $\wt W_\Delta=-W_\Gamma^{T}$. \ssq Le corollaire est alors une
cons\'equence directe du lemme \ref{type1} et du fait que l'on a $\wt
W_\Delta=\varepsilon(\Delta)\wt W_{\Delta^{\hbox{\sevenrm priv}}}$.
\end{preuve}\vspace{0,5cm}

Le cas d'une strate de type 1 se traite donc en posant
$\gamma_1=\alpha$, $\gamma_2=\gamma$, $\gamma_3=\beta$ et
$\gamma_j=\gamma$ pour $j\ge 4$. Donc avec les notations du lemme
le degr\'e de $\gamma_2$ vaut $2$.

Pour une strate de type 2; on \'echange donc les positions de
$\alpha$ et $\beta$ ce qui fournit un signe $(-1)^{k_1k_2}$. On
applique le lemme pour  $\gamma_1=\beta$ $\gamma_2=\gamma$,
$\gamma_3=\alpha$ et $\gamma_j=\gamma$ pour $j\ge 4$. On \'echange
a nouveau les position de $[\beta, \gamma]$ et $\alpha$ ce qui
produit  le signe $(-1)^{(k_2+1)k_1}$. Le produit des signes est
donc $(-1)^{k_1}$.  La somme des deux contributions fournit le
dernier terme dans le membre de droite de l'\'egalit\'e dans
l'\'enonc\'e de la proposition \ref{intermediaire}
\subsection{L'int\'egration sur une strate de type 3 a\'erien.}
Le corollaire \ref{cor} s'applique \'egalement aux strates de type 3
a\'erien en posant $\gamma_1=\gamma_2=\gamma$~: comme on a
$[\gamma,\,\gamma]=0$ la contribution de ces strates est nulle.
\subsection{L'int\'egration sur les strates de type 3 terrestre.}
\qquad Soit $Y_{3t}$ la r\'eunion des $Y\cap C_T$ o\`u $C_T$ est
une strate de type 3 terrestre (voir \S\ref{cylindre}). posons~:
$$W_\Gamma^{(3t)}=\int_{Y_{3t}}\omega_\Gamma.$$
On a alors~:
\begin{lemme}\label{type3t}
\begin{eqnarray}\nonumber
&\sum\limits_{n\ge 0}{\hbar^n\over n!}\sum\limits_{\Gamma\in
G_{n+2,m}} W_\Gamma^{(3t)}\Cal B_\Gamma
(\alpha\otimes\beta\otimes\gamma
        \otimes\cdots\otimes\gamma)(f_1\otimes\cdots\otimes f_m)=\\\nonumber
\hspace{-1cm}&=\sum\limits_{n\ge 0}{\hbar^n\over n!}
\sum\limits_{\Delta\in G_{n+2,m-1}} \wt W_\Delta\Cal
B_\Delta(\alpha\otimes\beta\otimes\gamma
        \otimes\cdots\otimes\gamma)
\bigl((f_1*f_2\otimes\cdots\otimes f_m)-\cdots\\
&+(-1)^{m-2}(f_1\otimes\cdots\otimes f_{m-1}*f_m)\bigr).
\end{eqnarray}
\end{lemme}
\begin{preuve}
Soit $\Gamma\in G_{n+2,m}$, soit $T$ une strate de type 3 terrestre
et soit $W_\Gamma^T$ l'int\'egrale de la forme $\omega_\Gamma$ sur
$C_T\cap Y$. Chaque sommet a\'erien du graphe interne est
\'etiquet\'e par le $2$-tenseur $\gamma$, donc il poss\`ede $2n_1$
fl\`eches, o\`u $n_1$ est le nombre de sommets a\'eriens de ce
graphe. Pour que $W_\Gamma^T$ soit non nul il faut qu'aucune
fl\`eche ne sorte du graphe interne, et il faut donc aussi que
$m_1=2$, o\`u $m_1$ est le nombre de points terrestres du nuage. Les
points terrestres du nuage sont donc $\{\overline i,\,\overline
{i+1}\}$ o\`u $i$ est un entier entre $1$ et $m-1$. On a donc ainsi
une partition de l'ensemble des strates du type 3 terrestre en $m-1$
groupes $E_i, i=1,\ldots ,m-1$. En sommant sur tous les graphes tels
que le graphe externe soit un $\Delta\in G_{n+2,m-1}$ donn\'e et
tels que le graphe interne s'obtient par \'eclatement du sommet
terrestre $\overline i$, on trouve donc~:
\begin{eqnarray}\nonumber
\hspace{-0,5cm}&\sum\limits_{\Gamma\in G_{n+p+2,m}\, , \,
\Gamma\to\Delta} {\hbar^p\over p!}{1\over|B_{n,p}|}
\sum\limits_{C_T\in E_i} W_\Gamma^T\Cal B_\Gamma
(\alpha\otimes\beta\otimes\gamma
        \otimes\cdots\otimes\gamma)(f_1\otimes\cdots\otimes f_m)=\\
\hspace{-1cm}&\pm \wt W_\Delta\Cal
B_\Delta(\alpha\otimes\beta\otimes\gamma
        \otimes\cdots\otimes\gamma)\bigl(f_1\otimes\cdots\otimes f_i*f_{i+1}\otimes\cdots\otimes f_m \bigr).
\end{eqnarray}
Comme les sommets $1$ et $2$ sont immuables il y a en effet
$|B_{n,p}|=(n+p)!/(n!p!)$ choix possibles d'ordre sur les sommets
a\'eriens de $\Gamma$ qui redonnent l'ordre du graphe interne et du
graphe externe (voir la discussion \`a la fin du \S\ 3.3). En
sommant sur tous les $\Delta$ dans $G_{n+2,m}$ \`a $n$ fix\'e, en
multipliant par $\hbar^n/n!$ et en sommant la s\'erie on obtient
alors~:
\begin{eqnarray}\nonumber
\sum\limits_{r\ge 0}{\hbar^r\over r!}\sum\limits_{\Gamma\in
G_{r+2,m}} \sum\limits_{C_T\in E_i} W_\Gamma^{T}\Cal B_\Gamma
(\alpha\otimes\beta\otimes\gamma
        \otimes\cdots\otimes\gamma)(f_1\otimes\cdots\otimes
        f_m)=\end{eqnarray}\begin{equation}
\pm\sum\limits_{n\ge 0}{\hbar^n\over n!} \sum\limits_{\Delta\in
G_{n+2,m-1}} \wt W_\Delta\Cal
B_\Delta(\alpha\otimes\beta\otimes\gamma
        \otimes\cdots\otimes\gamma)
\bigl(f_1\otimes f_2\otimes\cdots\otimes
f_i*f_{i+1}\otimes\cdots\otimes f_m).
\end{equation}
Pr\'ecisons maintenant le
signe $\pm$ en regardant les orientations~: si on d\'ecide
d'orienter $Y\cap C_T$ dans $C_T$ par sa normale sortante~:
$$\Omega_{C_T}=\dn\wedge \Omega_{Y\cap C_T},$$
on a alors~:
$$\Omega=\dn'\wedge\dn\wedge\Omega_{Y\cap C_T},$$
o\`u $\dn'$ est la normale sortante de la strate $C_T$ dans l'espace
de configuration $\overline{C_{n+2,m}}^+$. D'apr\`es \cite {AMM} \S\
I.2.2 la forme $\Omega_{C_T}$ est le produit ext\'erieur de la
``forme volume interne'' par la ``forme volume externe'' multipli\'e
par le signe $(-1)^{(i-1)2+(i-1)+2}=(-1)^{i-1}$. On obtient alors le
lemme (avec les signes altern\'es comme indiqu\'e) en sommant sur
tous les groupes de strates $E_i$.
\end{preuve}
On reconna\^\i t dans le membre de gauche un cobord de Hochschild
pour la multiplication d\'eform\'ee $*$, aux deux termes extr\^emes
pr\`es. Ces deux termes vont \^etre donn\'es par les strates de type
4~:
\subsection{L'int\'egration sur une strate de type 4.}
Soit $Y_4$ la r\'eunion des $Y\cap \Sigma$ o\`u $\Sigma$ est une
strate de type 4 (voir \S\ \ref{cylindre}). Posons~:
$$W_\Gamma^{(4)}=\int_{Y_{4}}\omega_\Gamma.$$
On a alors~:
\begin{lemme}\label{type4}
\begin{eqnarray}\nonumber
&\sum\limits_{n\ge 0}{\hbar^n\over n!}\sum\limits_{\Gamma\in
G_{n+2,m}} W_\Gamma^{(4)}\Cal B_\Gamma
(\alpha\otimes\beta\otimes\gamma
        \otimes\cdots\otimes\gamma)(f_1\otimes\cdots\otimes f_m)=\\\nonumber
&-\sum\limits_{n\ge 0}{\hbar^n\over n!}
\Bigl(\sum\limits_{\Delta\in G_{n+2,m-1}}f_1*\wt W_\Delta\Cal
B_\Delta(\alpha\otimes\beta\otimes\gamma
        \otimes\cdots\otimes\gamma)(f_2\otimes\cdots\otimes f_m)+\\
&(-1)^{m} \wt W_\Delta\Cal
B_\Delta(\alpha\otimes\beta\otimes\gamma
        \otimes\cdots\otimes\gamma)(f_1\otimes\cdots\otimes f_{m-1})*f_m\Bigr).
\end{eqnarray}
\end{lemme}
\begin{preuve}
Soit $\Gamma\in G_{n+2,m}$, soit $T$ une strate de type 4 et
soit $W_\Gamma^T$ l'int\'egrale de la forme $\omega_\Gamma$ sur
cette strate. Chaque sommet a\'erien du graphe externe est
\'etiquet\'e par le $2$-tenseur $\gamma$, donc il poss\`ede $2n_2$
fl\`eches, o\`u $n_2$ est le nombre de sommets a\'eriens de ce
graphe. Pour que $W_\Gamma^T$ soit non nul il faut qu'aucune
fl\`eche ne sorte du graphe interne, et il faut aussi que $m_2=2$,
o\`u $m_2$ est le nombre de points terrestres du graphe externe. Le
nuage contient donc exactement $m-1$ points terrestres. Les strates
qui v\'erifient ceci se r\'epartissent donc en deux classes, celles
qui laissent le premier sommet terrestre en-dehors du nuage, et
celles qui laissent le dernier. Un raisonnement similaire \`a celui
du paragraphe pr\'ec\'edent nous donne donc le lemme, aux signes \`a
pr\'eciser pr\`es. Les orientations nous donnent ici un signe
$(-1)^{(m-1)+(m-1)+1}=-1$ pour le premier terme, et $(-1)^{m-1}$
pour le deuxi\`eme terme.
\end{preuve}
\subsection{Fin de la d\'emonstration de la proposition \ref{intermediaire}}
La proposition \ref{intermediaire} d\'ecoule de l'addition des termes provenant
des quatre (ou plut\^ot cinq) types de strate (corollaire \ref{cor},
lemmes \ref{type3t} et \ref{type4}). L'addition des termes donn\'es par les
strates de type 3 et 4 donnent l'oppos\'e du cobord de Hochschild,
soit encore $[*,-]$.
\end{preuve}
\subsection{Fin de la d\'emonstration du th\'eor\`eme 1.2}
La formule de Stokes \ref{stokes} et les propositions
\ref{importante} et \ref{intermediaire} fournissent exactement le
r\'esultat cherch\'e \`a savoir le th\'eor\`eme 1.2 annonc\'e,
dont on peut pr\'esicer maintenant le cobord. On a donc :
\begin{theoreme}\label{final}
Soit $\alpha$ un $k_1$-champ de vecteurs et soit $\beta$ un
$k_2$-champ de vecteurs. Alors on a avec $m=k_1+k_2$~:
\begin{eqnarray}\nonumber
\Cal U'_{\hbar\gamma}(\alpha\cup\beta)- \Cal
U'_{\hbar\gamma}(\alpha)\cup \Cal U'_{\hbar\gamma}(\beta)=
\sum\limits_{n\ge 0}{\hbar^n\over n!}\!\sum\limits_{\Delta\in
G_{n+2,m-1}} \wt W_\Delta[*,\Cal
B_\Delta(\alpha\otimes\beta\otimes\gamma\otimes\cdots\otimes\gamma)]\end{eqnarray}
\begin{equation}
-\sum\limits_{n\ge 0}{\hbar^n\over n!}\sum\limits_{\Delta\in
G_{n+1,m}}\hskip -2mm \wt W_\Delta\Bigl(\Cal
B_\Delta([\alpha,\gamma]\otimes\beta\otimes
\gamma\otimes\cdots\otimes\gamma)+ (-1)^{k_1}\Cal
B_\Delta(\alpha\otimes[\beta,\gamma]\otimes\gamma\otimes\cdots\otimes\gamma)\Bigr).
\end{equation}
\end{theoreme}
Si on a
$[\alpha,\,\gamma]=[\beta,\,\gamma]=0$ le dernier terme est nul, et
$U'_{\hbar\gamma}(\alpha\cup\beta)-\Cal U'_{\hbar\gamma}(\alpha)\cup
\Cal U'_{\hbar\gamma}(\beta)$ est bien donn\'e par un cobord de
Hochschild pour la multiplication d\'eform\'ee $*$.
%
\bibliographystyle{plain}
\bibliography{mt}
\end{document}